\documentclass{article}

\usepackage[utf8]{inputenc}
\usepackage[T1]{fontenc}
\usepackage[a4paper,margin=1.75cm]{geometry}
\usepackage{hyperref}
\hypersetup{bookmarksopen=true,colorlinks=true,allcolors={blue},}
\usepackage{authblk}
\usepackage{amssymb}
\usepackage{amsthm}
\usepackage{amsfonts}
\usepackage{graphicx}
\usepackage{float}
\usepackage{subcaption}
\usepackage{tabu}
\usepackage{multirow}
\usepackage{mathtools}
\usepackage{cleveref}

\theoremstyle{plain}% default, other are: remark, definition
\newtheorem{theorem}{Theorem}[section]

\newtheorem{remark}[theorem]{Remark}

\numberwithin{equation}{section}
\usepackage{pifont}
\newcommand{\tick}{\ding{51}}
\newcommand{\cross}{\ding{55}}
\DeclareMathOperator\real{Re}

\renewcommand\Re{\real}

\newcommand{\assign}{\mathrel{\mathop:}=}

\usepackage{pgfplotstable,booktabs}
\pgfplotsset{compat=newest}

\usepackage[
    sorting=nyt,
    minbibnames=5,
    maxbibnames=5,
    giveninits=true,
    url=true,
    isbn=false,
    eprint=true,
    doi=true
]{biblatex}

\title{Achieving wavenumber robustness\\
in domain decomposition for heterogeneous Helmholtz equation:\\
an overview of spectral coarse spaces}

\author[,1]{Victorita Dolean\thanks{v.dolean.maini@tue.nl}}
\author[,2]{Mark Fry\thanks{mark.fry@strath.ac.uk}}
\author[,2]{Matthias Langer\thanks{m.langer@strath.ac.uk}}
\author[,3]{Emile Parolin\thanks{emile.parolin@inria.fr}}
\author[,3]{Pierre-Henri Tournier\thanks{pierre-henri.tournier@sorbonne-universite.fr}}

\affil[1]{Eindhoven University of Technology, PO Box 513, Eindhoven, 5600 MB, The Netherlands}
\affil[2]{University of Strathclyde, 26 Richmond Street, Glasgow, G1 1XH, United Kingdom}
\affil[3]{Sorbonne Université, Université Paris Cité, CNRS, INRIA, Laboratoire Jacques-Louis Lions, LJLL, EPC ALPINES, 4 place Jussieu, F-75005 Paris, France}

\date{}

\begin{document}

\maketitle

\begin{abstract}
Solving time-harmonic wave propagation problems in the frequency domain within heterogeneous media poses significant mathematical and computational challenges, particularly in the high-frequency regime. Among the available numerical approaches, domain decomposition methods are widely regarded as effective due to their suitability for parallel computing and their capacity to maintain robustness with respect to physical parameters, such as the wavenumber. These methods can achieve near-constant time-to-solution as the wavenumber increases, though often at the expense of a computationally intensive coarse correction step. This work focuses on identifying the best algorithms and numerical strategies for benchmark problems modelled by the Helmholtz equation. Specifically, we examine and compare several coarse spaces which are part of different families, e.g.\ GenEO (Generalised Eigenvalue Overlap) type coarse spaces and harmonic coarse spaces, that underpin two-level domain decomposition methods. By leveraging spectral information and multiscale approaches, we aim to provide a comprehensive overview of the strengths and weaknesses of these methods. Numerical experiments demonstrate that the effectiveness of these coarse spaces depends on the specific problem and numerical configuration, highlighting the trade-offs between computational cost, robustness, and practical applicability.
\end{abstract}



\section{Introduction}
The main focus of this work is tackling the computational complexity inherent in frequency-domain simulations of wave propagation and scattering for heterogeneous problems. These challenges are prevalent in various engineering fields, including acoustics, electromagnetic analysis, and seismic imaging. The finite element discretisation of frequency-domain wave models typically produces large-scale, indefinite, and poorly conditioned linear systems, which are especially difficult to handle with standard solvers. This difficulty escalates with high-frequency scenarios and intricate heterogeneities. Ensuring precision requires an increase in grid resolution at least proportional to frequency, making the resulting systems infeasible for direct solution methods in high-frequency cases. Consequently, the need arises for carefully tailored iterative techniques. This work focuses on two-level domain decomposition strategies aimed at efficiently solving these linear systems in parallel.

{The wave behaviour is described by a heterogeneous Helmholtz equation posed on a computational domain 
\(\Omega \subset \mathbb{R}^d\) (\(d = 2,3\)), where the complex-valued field 
\(u(\boldsymbol{x})\) satisfies  
\begin{alignat}{3}\label{HelmholtzSystemUnified}
	- \nabla \cdot \big( a(\boldsymbol{x}) \nabla u \big) - \omega^2\, m(\boldsymbol{x})\, u
	&= f (\boldsymbol{x})
	\hspace*{10ex} &&\text{in } \Omega, \\
	\mathcal{C}(u)
	&= 0
	\hspace*{10ex} &&\text{on } \partial \Omega,
\end{alignat}
with \(\mathcal{C}\) representing the imposed boundary conditions, 
\(a(\boldsymbol{x}) > 0\) a spatially varying coefficient associated with material stiffness 
(or inverse density), \(m(\boldsymbol{x}) > 0\) a spatially varying ``mass'' coefficient related 
to the squared slowness, \(\omega > 0\) the angular frequency,
$c(\boldsymbol{x})$ the wave speed and \(f(\boldsymbol{x})\) the source term.
}

{The local wavenumber is given by
\[
k(\boldsymbol{x}) = \omega \sqrt{\frac{m(\boldsymbol{x})}{a(\boldsymbol{x})}},
\]
which reduces to the familiar relation \(k(\boldsymbol{x}) = \omega/c(\boldsymbol{x})\) when 
\(a \equiv 1\) and \(m(\boldsymbol{x}) = c^{-2}(\boldsymbol{x})\).  
In heterogeneous media, spatial variations in \(a\) and \(m\) lead to additional variations in 
\(k(\boldsymbol{x})\), making the high-frequency regime (large \(k\)) particularly 
challenging for numerical solvers.}

The challenges of developing effective solvers for the Helmholtz equation are well documented in reviews such as \cite{Ernst:2012:NAM, Gander:2018:SIREV}, which highlight the difficulty in extending state-of-the-art methods for symmetric positive definite problems to handle the indefinite or non-self-adjoint nature of the Helmholtz problem. For large-scale problems---where accurate discretisation for high wavenumbers is required---domain decomposition methods emerge as a natural approach \cite{Dolean:15:DDM}. Despite recent advances, both theoretical \cite{Graham:2017:RRD, Graham:2018:DDI, Gong:2021:DDP} and numerical \cite{Dolean:2020:LFD, Dolean:2020:IFD, Tournier:2022:3FD, Operto:2023:I3F}, we are still at the stage where \textit{no single method has consistently outperformed others in addressing the solution of the discretised Helmholtz equation in high frequency regime}. 

Domain decomposition methods are widely used for solving large systems arising from PDE discretisation, offering robust strategies for various problems. However, traditional methods often fail or even diverge when applied to wave propagation problems. One critical aspect requires careful attention: the coarse space that captures global behaviour and facilitates distant inter-domain information transfer while dealing with the non-locality of the wave phenomena. This work focuses on overlapping Schwarz methods which are robust against the physical coefficients (wavenumber and heterogeneity) of the problem and scalable as the number of subdomains increases. Addressing these challenges requires coarse information that is both computationally efficient and globally accessible to all subdomains. This study emphasises once again the importance of sophisticated coarse spaces specially tailored to the nature of the problem. 

One of the prominent methods to tackle complex problems with large contrasts in the physical coefficients is the GenEO (Generalised Eigenproblems in the Overlap) coarse space which has proven effective for solving highly heterogeneous elliptic problems \cite{Spillane:2014:ARC, Haferssas:2017:ADS}. However, for the Helmholtz equation, selecting an appropriate coarse space is particularly challenging. Due to its indefinite or non-self-adjoint nature, enlarging the coarse space does not necessarily enhance performance \cite{Fish:2000:GBT}. Designing coarse spaces for the Helmholtz problem also involves minimising dependence on the wavenumber $k$. Using plane waves as a coarse space basis is a natural approach for capturing global behaviour, but their suitability for heterogeneous media remains uncertain. Plane waves were first applied in multigrid methods \cite{Brandt:1997:WRM} and later extended to domain decomposition methods, such as FETI(-DP)-H \cite{Farhat:2000:ATL, Farhat:2005:FDP}, though primarily for homogeneous problems. 

Even though coarse space information must be globally available, coarse spaces can be constructed locally using basis vectors derived from solutions to local eigenvalue problems. Spectral coarse spaces, such as the DtN (Dirichlet-to-Neumann) coarse space, illustrate this approach for the Helmholtz equation \cite{Conen:2014:ACS}. These eigenproblems are formulated on subdomain interfaces using a DtN map, building on methods developed for elliptic problems \cite{Nataf:2011:ACS, Dolean:2012:ATL}. Due to the difficulty in obtaining theoretical results for these kinds of problems, the first assessment  of two different yet related spectral coarse spaces---the DtN method and a GenEO-type approach tailored to the Helmholtz problem---was done in \cite{Bootland:2022:CDG} for simple test cases. At the same time in \cite{Bootland:2021:ACS} the authors have synthetised recent advancements in enhancing domain decomposition methods for heterogeneous media and implemented these approaches using a common software framework, FreeFEM. Practical implementation details and trade-offs were thoroughly discussed there, and extensive numerical experiments on 2D and 3D benchmark problems were presented to evaluate and compare the methods in various settings. Finally, those results provide guidance on the suitability of different methods for specific scenarios, offering insights for future applications.

More recently, a new family of coarse spaces has emerged, inspired by the multiscale spectral generalised finite element method (MS-GFEM). In \cite{Ma:2024:TLR} the authors introduce global coarse solves derived from MS-GFEM spaces constructed through local eigenproblems and prove that the preconditioner, when used with GMRES, achieves convergence at a rate tied to the MS-GFEM approximation error, and since MS-GFEM exhibits exponential convergence due to oversampling this will ensure efficient performance with a small coarse space.  This work provides a robust and scalable strategy for solving heterogeneous Helmholtz problems, particularly in high-frequency regimes. Akin to this philosophy \cite{Hu:2024:ANC} relies on coarse spaces based on generalised eigenvalue problems on harmonic spaces composed of Helmholtz extensions.
Besides, \cite{Galkowski:2025:CTT} provides a general theoretical framework of
two-level methods for Helmholtz.
Another approach has been described in~\cite{Nataf2024}, and provides a
systematic way of constructing coarse spaces for possibly non-Hermitian and
indefinite matrices hence applicable to the linear systems arising from the
discretisation of the Helmholtz equation.
Finally, we mention~\cite{Daas2024} which provides a fully algebraic coarse
space construction that has been applied to indefinite or non-symmetric
matrices, with available analysis only for SPD matrices.

Although the previous methods are theoretically sound, they may face several limitations that may impact their practical applicability. The convergence theory relies on assumptions, which may not hold in real-world scenarios. While methods are introduced in a general setting, their scalability for large-scale 3D problems and heterogeneous media remains unproven. Additionally, the practical utility of coarse spaces, including both problem-adapted and non-adapted basis functions, is uncertain, as the results may not translate into computationally efficient implementations. The focus on theoretical rigor over implementation leaves challenges like the cost of coarse space construction and computational bottlenecks unaddressed. Furthermore, the preconditioner’s performance may still be heavily influenced by problem-specific parameters, raising concerns about robustness in practical settings. 

{Our main contributions are:
\begin{itemize}
\item We provide the first unified and reproducible comparison of all major spectral coarse space families for the heterogeneous, high‑frequency Helmholtz equation---implemented within a common FreeFEM framework.
\item We carry out large‑scale numerical experiments on both simple test cases and challenging 2D/3D benchmarks, revealing new insights into robustness, scalability, and computational cost.
\item We provide practical guidance on when and how each coarse space is most effective, offering actionable recommendations for real‑world Helmholtz solvers.
\end{itemize}}

{Building on the earlier comparative study of coarse spaces for the Helmholtz equation in the high-frequency regime \cite{Bootland:2021:ACS}, this work substantially extends the scope and depth of the analysis. In particular, we incorporate the most recent advances in coarse space design---most notably the emerging family of Helmholtz-harmonic and extended harmonic coarse spaces---and examine them alongside established DtN and GenEO-type approaches. The numerical evaluation goes far beyond previous work: we explore performance across a wide spectrum of 2D and 3D benchmark problems, from idealised homogeneous cases to challenging heterogeneous, high-contrast, and realistic large-scale geometries, pushing the limits of current two-level methods. This leads to new insights into the trade-offs between robustness, coarse space size, and computational cost, and clarifies the conditions under which each approach is most effective.}

{The paper is organised as follows. \Cref{sec:HelmholtzProblem} introduces the heterogeneous Helmholtz problem, its boundary value formulation, and finite element discretisation. \Cref{sec:DomainDecomposition} outlines the foundations of domain decomposition methods and motivates the introduction of a second-level (coarse space) correction via deflation. \Cref{sec:coarse} presents the three main families of coarse spaces---DtN, GenEO-type, and harmonic---together with their recent variants. \Cref{sec:ComparativeNumericalStudies} provides an extensive set of numerical experiments, covering both controlled test cases and realistic large-scale applications. Conclusions and perspectives are discussed in \Cref{sec:Conclusions}.}

\section{Discretisation of the heterogeneous Helmholtz problem}
\label{sec:HelmholtzProblem}
We focus here on the interior heterogeneous Helmholtz problem 
\eqref{HelmholtzSystemUnified} with suitable boundary conditions. 
In practical settings, {the domain is truncated to a bounded computational domain $\Omega$},
and the far-field Sommerfeld radiation condition must be approximated on the artificial boundary of $\Omega$ to model wave behaviour appropriately.  
A widely used approximation is the Robin (impedance) condition, though  alternatives include absorbing boundary conditions (ABCs) 
\cite{Zarmi:2013:AGA} or perfectly matched layers (PML) 
\cite{Harari:2000:AAN,Beriot:2021:AAP}.  
We consider Dirichlet conditions on $\Gamma_{D} \subset \partial\Omega$ 
and Robin conditions on $\Gamma_{R} = \partial\Omega \setminus \Gamma_{D}$, 
solving the boundary value problem:
{
\begin{subequations}
\label{HelmholtzSpecificSystemUnified}
\begin{alignat}{2}
\label{HelmholtzSpecificEquationUnified}
- \nabla \cdot \big( a(\boldsymbol{x}) \nabla u \big) 
- \omega^2 m(\boldsymbol{x})\, u &= f(\boldsymbol{x}) \hspace*{10ex}
& & \text{in } \Omega,\\
\label{HelmholtzSpecificDirichletBCUnified}
u &= u_{\Gamma_{D}} & & \text{on } \Gamma_{D},\\
\label{HelmholtzSpecificRobinBCUnified}
a(\boldsymbol{x})\,\frac{\partial u}{\partial \boldsymbol{n}} 
+ i\,\omega\,Z(\boldsymbol{x})\,u &= 0 & & \text{on } \Gamma_{R},
\end{alignat}
\end{subequations}
where {$\boldsymbol{n}$ is the outward normal,}
and $Z(\boldsymbol{x})$ is an impedance coefficient typically chosen as 
$Z(\boldsymbol{x}) = \sqrt{a(\boldsymbol{x})\,m(\boldsymbol{x})}$ 
to approximate the Sommerfeld radiation condition.  
The problem is well posed if $\Gamma_{R} \neq \emptyset$ but may become 
ill-posed for certain parameter choices when $\Gamma_{R} = \emptyset$. The local wavenumber $k(\boldsymbol{x})$ can vary within $\Omega$ and may have discontinuities across 
material interfaces, allowing $k \in L^{\infty}(\Omega)$. } 

For \(\Omega' \subset \Omega\), the trial and test spaces are defined as
\begin{equation}
\label{eq:vomega}
V(\Omega') = \{ u \in H^{1}(\Omega') : u = u_{\Gamma_{D}} 
\text{ on } \Gamma_{D} \cap \partial \Omega' \}, \qquad
V_{0}(\Omega') = \{ u \in H^{1}(\Omega') : u = 0 
\text{ on } \Gamma_{D} \cap \partial \Omega' \}.
\end{equation}
The weak formulation seeks $u \in V(\Omega)$ satisfying:
\begin{equation}
\label{WeakFormUnified}
a_{\Omega}(u,v) = F(v), \quad \forall v \in V_{0}(\Omega),
\end{equation}
where, for any \(\Omega' \subset \Omega\), we define:
\begin{align}
\label{WeakFormTermsUnified}
a_{\Omega'}(u,v) &= 
\int_{\Omega'} \big( a(\boldsymbol{x}) \nabla u \cdot \nabla \bar{v} 
- \omega^2 m(\boldsymbol{x})\,u\,\bar{v} \big) \, \mathrm{d}\boldsymbol{x}  + \int_{\Gamma_{R} \cap \partial\Omega'} i\,\omega\,Z(\boldsymbol{x})\,u\,\bar{v} \, \mathrm{d}s, \qquad
F(v) = \int_{\Omega} f\,\bar{v} \, \mathrm{d}\boldsymbol{x}.
\end{align}

Using a simplicial mesh $\mathcal{T}^{h}$ of $\Omega$ with element diameter
$h$, we use finite elements to discretise \eqref{WeakFormUnified}, 
as outlined in \cite{Bootland:2021:ACS}.  
A Lagrange finite element approximation reduces the problem to solving the
complex linear system:
\begin{equation}
\label{LinearSystemUnified}
A\mathbf{u} = \mathbf{f},
\end{equation}
where $A \in \mathbb{C}^{n \times n}$ and $\mathbf{f} \in \mathbb{C}^{n}$ are respectively derived 
from the sesquilinear and antilinear forms, and \(n\) denotes the
total number of degrees of freedom in \(\Omega\).  
Accurate approximation requires fine discretisation, particularly as 
$k(\boldsymbol{x})$ increases, due to the pollution effect \cite{Babuska:1997:IPE}.  
For piecewise linear (\(\mathbb{P}_{1}\)) elements, maintaining accuracy requires 
$k_{\max}^3 h^2$ to remain bounded, where 
$k_{\max} = \sup_{\mathbf{x}\in\Omega} k(\boldsymbol{x})$, 
implying $h = {\cal O} (k_{\max}^{-3/2})$.  
Higher-order elements relax this constraint, but interpolation properties 
degrade for very high orders.  
An alternative common strategy is to fix the number of points per wavelength, 
$n_\text{ppwl}$, with $h = \mathcal{O}(k_{\max}^{-1})$, 
where the wavelength is $\lambda(\boldsymbol{x}) = 2 \pi / k(\boldsymbol{x})$.  
Practitioners often use 5 or 10 points per wavelength.

\section{Domain decomposition methods}
\label{sec:DomainDecomposition}
{To solve the discrete heterogeneous Helmholtz problem 
\eqref{LinearSystemUnified}, we employ GMRES with acceleration provided by a two-level overlapping domain decomposition preconditioner.  The core one-level approach used is the optimised restricted additive Schwarz 
(ORAS) method \cite{Dolean:15:DDM}.  To formulate the domain decomposition preconditioner, we begin by dividing 
the domain $\Omega$ into a set of \(N\) non-overlapping subdomains 
$\{ \Omega_s' \}_{s=1}^N$, which align with the finite element mesh  $\mathcal{T}^h$.  
Overlapping subdomains $\{ \Omega_s \}_{s=1}^N$ are then created by expanding  each $\Omega_s'$ with neighbouring mesh elements, defined as:
\[
\Omega_s = \mathrm{Int}\left(
\bigcup_{\mathrm{supp}(\phi_j) \cap \Omega_s' \neq \emptyset} 
\mathrm{supp}(\phi_j)
\right),
\]
where $\mathrm{Int}(\cdot)$ and $\mathrm{supp}(\cdot)$ represent the interior 
and support of a function, respectively, and $\phi_j$ are the nodal basis 
functions of the finite element space.  
Additional layers can be recursively added to increase the overlap as needed.}

{With the overlapping decomposition, we define the preconditioner components.  
The restriction operator $R_s \in \mathbb{R}^{n_s \times n}$ 
restricts functions to $\Omega_s$ where \(n_{s}\) denotes the local number of
degrees of freedom in \(\Omega_{s}\), while the transpose $R_s^T$ extends
by zero outside $\Omega_s$.  
A partition of unity is represented by diagonal matrices 
$D_s \in \mathbb{R}^{n_s \times n_s}$, ensuring
\[
\sum_{s=1}^N R_s^T D_s R_s = I.
\]}
{ORAS also requires solving local heterogeneous Helmholtz–Robin problems:
\begin{subequations}
\label{ORASLocalSystemUnified}
\begin{alignat}{2}
- \nabla \cdot \big( a(\boldsymbol{x}) \nabla w_s \big) 
- \omega^2 m(\boldsymbol{x})\, w_s &= f \hspace*{10ex}
& & \text{in } \Omega_s, \\
a(\boldsymbol{x})\,\frac{\partial w_s}{\partial \boldsymbol{n}_s} 
+ i\,\omega\,Z(\boldsymbol{x})\,w_s &= 0 
& & \text{on } \partial\Omega_s \setminus \partial\Omega, \\
\mathcal{C}(w_s) &= 0 
& & \text{on } \partial\Omega_s \cap \partial\Omega,
\end{alignat}
\end{subequations}
where $\mathcal{C}$ applies the boundary conditions 
\eqref{HelmholtzSpecificDirichletBCUnified}--%
\eqref{HelmholtzSpecificRobinBCUnified}, 
and $Z(\boldsymbol{x})$ is the local impedance, typically chosen as 
$Z(\boldsymbol{x}) = \sqrt{a(\boldsymbol{x})\,m(\boldsymbol{x})}$.}  The finite element discretisation of \eqref{ORASLocalSystemUnified} yields a 
local stiffness matrix $\widehat{A}_s \in \mathbb{C}^{n_s \times n_s}$.  
The ORAS preconditioner is then assembled as:
\[
M_{\text{ORAS}}^{-1} = 
\sum_{s=1}^N R_s^T D_s \widehat{A}_s^{-1} R_s,
\]
where the local solves $\widehat{A}_s^{-1}$ are performed in parallel.

{To ensure robustness and scalability, particularly for indefinite problems 
like the Helmholtz equation, the ORAS method incorporates a coarse space 
in a two-level framework.  
The coarse space, represented by a matrix $Z$ with linearly independent 
column vectors, plays a key role in addressing scalability.  
Using deflation, a coarse operator $E = Z^\dagger A Z$ and correction 
operator $Q = Z E^{-1} Z^\dagger$ are defined, yielding the two-level 
ORAS method:
\[
M_{\text{ORAS,2}}^{-1} = 
M_{\text{ORAS}}^{-1}(I - AQ) + Q.
\]
The choice of coarse space is critical, particularly for indefinite problems 
where its addition may not always improve performance 
\cite{Fish:2000:GBT}.  
We next discuss the selection of coarse spaces.}

\section{Spectral coarse spaces}
\label{sec:coarse}
In this study, we investigate spectral coarse spaces for the discrete Helmholtz problem \eqref{LinearSystemUnified}, which are constructed by solving local eigenvalue problems on subdomains and assembling the resulting eigenfunctions into a global coarse space. While a broad family of such spectral coarse spaces exists, we focus on two representative classes that follow a similar construction philosophy: {\bf GenEO-type spaces} and {\bf harmonic-type spaces}.

We begin by reviewing the DtN coarse space introduced in \cite{Conen:2014:ACS}, then present the $\Delta$-GenEO \cite{Bootland:2022:OSM, Dolean:2024:ITE} and $H_k$-GenEO spaces \cite{Dolean:2024:SPW_updated}, which are inspired by the GenEO methodology \cite{Spillane:2014:ARC}. In parallel, we consider recently developed Helmholtz-harmonic coarse spaces from \cite{Hu:2024:ANC, Ma:2024:TLR, Nataf2024}. Finally, we establish connections between these approaches.

Although comparative analyses were conducted in \cite{Bootland:2021:ACS, Bootland:2022:CDG}, our study builds upon and extends them by incorporating the latest developments in the field, along with a detailed complexity analysis on both model and benchmark problems. Note that in light of the same recent developments the DtN coarse space can be both considered as an ancestor of GenEO coarse spaces but also, by construction, falls into the category of harmonic coarse spaces.

\textbf{Notation}: Given a variational problem with system matrix $B$, we denote by $B_s$ the local Dirichlet matrix on subdomain $\Omega_s$. If Robin boundary conditions are imposed on internal interfaces, the corresponding local matrix is denoted by $\widehat{B}_s$, whereas for Neumann boundary conditions we write $\widetilde{B}_s$.

\subsection{The DtN Coarse space}

We begin by revisiting the Dirichlet-to-Neumann (DtN) coarse space, originally introduced in \cite{Nataf:2011:ACS} for elliptic problems and extended to the Helmholtz equation in \cite{Conen:2014:ACS}. This method constructs the coarse space from the low-frequency eigenmodes of a boundary operator defined via local DtN maps, which are then extended harmonically inside each subdomain. From a broader perspective, the DtN approach can be viewed as an early instance of harmonic-type coarse spaces, and also as a precursor to GenEO-type constructions that build on similar principles of localised spectral analysis.

{Let $\Gamma_s = \partial \Omega_s \setminus \partial \Omega$ denote the 
interface between subdomain $\Omega_s$ and its neighbours.  
Given Dirichlet data $v_{\Gamma_s}$ on $\Gamma_s$, 
\emph{the local heterogeneous Helmholtz extension} is the solution $v$ of 
the boundary value problem:
\begin{subequations}
\label{HelmholtzExtensionUnified}
\begin{alignat}{2}
- \nabla \cdot \big( a(\boldsymbol{x}) \nabla v \big) 
- \omega^2 m(\boldsymbol{x})\, v &= 0 \hspace*{10ex}
& & \text{in } \Omega_s, \\
v &= v_{\Gamma_s} & & \text{on } \Gamma_s, \\
\mathcal{C}(v) &= 0 & & \text{on } \partial\Omega_s \cap \partial\Omega,
\end{alignat}
\end{subequations}
where $\mathcal{C}(v) = 0$ enforces the original problem boundary conditions. The DtN map then takes $v_{\Gamma_s}$ to its corresponding Neumann data on 
$\Gamma_s$:
\begin{equation}
\label{DtNMapUnified}
\mathrm{DtN}_{\Omega_s}(v_{\Gamma_s}) = 
\left. a(\boldsymbol{x})\,\frac{\partial v}{\partial n} \right\rvert_{\Gamma_s}.
\end{equation}}
The DtN eigenproblem for subdomain $\Omega_s$ is
\begin{equation}
\label{DtNEigenproblemUnified}
\mathrm{DtN}_{\Omega_s}(u_{\Gamma_s}) = \lambda\,u_{\Gamma_s},
\end{equation}
where $u_{\Gamma_s}$ are eigenfunctions with eigenvalues 
$\lambda \in \mathbb{C}$.  
The corresponding coarse space is constructed by taking the heterogeneous 
Helmholtz extension of $u_{\Gamma_s}$ in $\Omega_s$ and extending it by zero 
across $\Omega$ via the partition of unity \cite{Conen:2014:ACS}.

For the discrete setting, we introduce local Neumann matrices $\widetilde{A}_s$, where $\mathcal{C} = 0$ on $\partial\Omega_s\cap\partial\Omega$. Let $\Gamma_s$ and $I_s$ denote boundary and interior indices, respectively. Defining the mass matrix on the interface as $\label{DtNMassMatrix} M_{\Gamma_s} = \left(\int_{\Gamma_s} \phi_j \phi_i \right)_{i,j \in \Gamma_s}$ the discrete DtN eigenproblem is given by 
\begin{align} \label{DiscreteDtNEigenproblem} 
\left(\widetilde{A}_{\Gamma_s,\Gamma_s} -A_{\Gamma_s,I_s}A_{I_s,I_s}^{-1}A_{I_s,\Gamma_s}\right) \mathbf{u}_{\Gamma_s} = \lambda M_{\Gamma_s} \mathbf{u}_{\Gamma_s},
\end{align}
{where $A_{I_s,\Gamma_s}, A_{\Gamma_s,I_s}$ and $A_{I_s,I_s}, A_{\Gamma_s,\Gamma_s}$ are obvious notations for sub-blocks of matrix $A$ corresponding to the set of indices in the interior $I_s$ and on the boundary $\Gamma_s$.} 
Because of the Robin boundary condition on \(\Gamma_{R}\), the
left-hand-side matrix in~\eqref{DiscreteDtNEigenproblem} is a priori
non-Hermitian.
The Helmholtz extension to interior degrees of freedom is obtained via
$\mathbf{u}_{I_s} = - A_{I_s,I_s}^{-1}A_{I_s,\Gamma_s} \mathbf{u}_{\Gamma_s}$,
and the global coarse basis in \(Z\) is formed as $R_s^T D_s \mathbf{u}_s$, where
$\mathbf{u}_s$ is the full local vector computed from
\(\mathbf{u}_{\Gamma_{s}}\) and \(\mathbf{u}_{I_{s}}\) \cite{Conen:2014:ACS}.
We see that the discrete counterpart of \Cref{DtNEigenproblemUnified} is
\Cref{DiscreteDtNEigenproblem} consisting in a Schur complement based
generalised eigenvalue problem.

A key aspect is the selection of eigenvectors for the coarse space. In \cite{Conen:2014:ACS} a criterion has been identified empirically consisting in retaining the eigenvalues with the smallest real parts, imposing the threshold \begin{equation} \label{Threshold} 
\Re(\lambda) < \eta_{\text{max}}. 
\end{equation} 
While $\eta_{\text{max}} = k_s$, where $k_s\assign \sup_{\boldsymbol{x}\in\Omega_s}k(\boldsymbol{x})$, was initially proposed, recent work \cite{Bootland:2019:ODN, Bootland:2022:CDG} suggests that $\eta_{\text{max}} = k_s^{4/3}$ can improve robustness with respect to the wavenumber. Also in \cite{Bootland:2022:CDG} this threshold was systematically studied and empirical $k$-dependent bounds were derived for the sizes of the corresponding coarse spaces.

\begin{remark}
Unlike for Laplace's equation {with self-adjoint boundary conditions} where generalised eigenvalue problems are self-adjoint and positive definite, the DtN eigenproblems for Helmholtz are in general non-self-adjoint. As a consequence, there are currently no theoretical convergence guarantees for the DtN coarse space. Nevertheless, recent studies---including \cite{Bootland:2021:ACS} and \cite{Bootland:2022:CDG}---demonstrate that this approach performs remarkably well in practice, and its complexity can be quantified empirically. A significant advantage of this construction is that its coarse space dimension is inherently limited by the number of degrees of freedom on the subdomain interface~$\Gamma_s$, as the eigenfunctions are defined on~$\Gamma_s$ and only extended into the interior. This allows for global modes with full support to be generated from a relatively compact spectral basis.
\end{remark}

\subsection{GenEO-type coarse spaces: $\Delta$-GenEO and $H$-GenEO}
The Generalised Eigenproblems in the Overlap (GenEO) coarse space, introduced in \cite{Spillane:2014:ARC}, offers a robust approach for symmetric positive definite problems, even in heterogeneous settings. 
Note that originally, GenEO was formulated variationally under the assumption that the bilinear form $a_{\Omega}(\cdot,\cdot)$ is symmetric and coercive, assumption not satisfied here. Given a subdomain $\Omega_s$ and its overlapping region $\Omega_s^\circ$, the local eigenproblem is defined as 
\begin{align} \label{GenEOVariationalEigenproblem} 
a_{\Omega_s}(u,v) = \lambda a_{\Omega_s^\circ}(\Xi_s(u),\Xi_s(v)) \quad \forall v \in V(\Omega_s), 
\end{align} 
where $\Xi_s$ is the partition of unity operator. Here, $a_D(\cdot,\cdot)$ represents the variational problem restricted to domain $D$ with problem boundary conditions on $\partial\Omega$ and natural conditions on $\partial D \setminus \partial\Omega$. The smallest eigenvalue $\lambda$ not included in the GenEO space determines the preconditioned operator’s condition number, making thresholding on $\lambda < \lambda_{\text{max}}$ a natural choice. Alternatively, a fixed number of eigenvectors per subdomain can be selected for efficiency. The restriction to $\Omega_s^\circ$ is not necessary, and a more practical variant replaces it with $\Omega_s$, simplifying implementation \cite{Dolean:15:DDM}. This leads to the discrete eigenproblem \begin{align} \label{DiscreteGenEOEigenproblem} \widetilde{A}_s \mathbf{u} = \lambda D_s A_s D_s \mathbf{u}, \end{align} where $A_s = R_s A R_s^T$ is the local Dirichlet matrix. The coarse space basis in \(Z\) is then constructed as $R_s^T D_s \mathbf{u}$. This formulation will serve as the foundation for a Helmholtz-specific extension.

\subsubsection{$\Delta$-GenEO coarse space}
Applying GenEO to the Helmholtz problem is challenging due to the loss of self-adjointness and definiteness. While theoretical advancements have been made \cite{Bootland:2022:OSM,Bootland:2022:GCS}, a fully justified spectral coarse space for Helmholtz remains elusive. One workaround, $\Delta$-GenEO, replaces the Helmholtz operator with the Laplacian (setting $\omega=0$ in \Cref{HelmholtzSpecificSystemUnified} and \Cref{WeakFormTermsUnified}), yielding the eigenproblem 
\begin{align} \label{DiscreteLaplaceGenEOEigenproblem} 
\widetilde{L}_s \mathbf{u} = \lambda D_s L_s D_s \mathbf{u}, 
\end{align} 
where $L_s$ and $\widetilde{L}_s$ are the local Laplace Dirichlet and Neumann matrices, respectively. This provides a well-posed problem with real non-negative eigenvalues, but its effectiveness deteriorates as $k$ increases since Laplace solutions diverge significantly from those of Helmholtz. For this coarse space theoretical bounds have been rigorously established in \cite{Bootland:2022:OSM} and further improved in \cite{Dolean:2024:ITE}. These bounds show, as expected, that this coarse space cannot perform very effectively at high frequencies and the size of the coarse space increases dramatically with the frequency.

\subsubsection{The $H_k$-GenEO coarse space}

To construct a more suitable spectral coarse space for the heterogeneous 
Helmholtz problem, it is essential to incorporate the Helmholtz operator 
itself.  
However, directly applying \Cref{DiscreteGenEOEigenproblem} with the 
heterogeneous operator leads to a non-self-adjoint problem with complex 
eigenvalues, which can cause solver instability.  
Instead, we propose a hybrid approach that links the indefinite heterogeneous 
Helmholtz operator to a positive-definite surrogate operator.  The local eigenproblem is given by
\begin{equation}
\label{DiscreteH-GenEOEigenproblemUnified}
\widetilde{A}_s \mathbf{u} 
= \lambda\, D_s\, P_s\, D_s \mathbf{u},
\end{equation}
where:
{
\begin{itemize}
\item $\widetilde{A}_s$ is the local Neumann matrix for the heterogeneous 
Helmholtz operator
\[
- \nabla \cdot \big( a(\boldsymbol{x}) \nabla u \big) 
- \omega^2 m(\boldsymbol{x})\, u,
\]
\item $D_s$ is the local partition-of-unity weighting matrix,
\item $P_s$ is a symmetric positive-definite (SPD) matrix chosen to be 
``close'' to $A_s$ in a suitable sense.
\end{itemize}
A natural choice for $P_s$ is $L_s$, the discretisation of the Laplace 
operator $- \Delta$ with heterogeneous coefficient $a(\boldsymbol{x})$.  
In practice, many SPD choices perform well, but another particularly 
convenient choice---allowing for simpler theoretical estimates---is
\(P_s = B_s\),
where $B_s$ is the discretisation of the \emph{positive} heterogeneous Helmholtz operator
\begin{equation}
\label{eq:posHelm}
- \nabla \cdot \big( a(\boldsymbol{x}) \nabla u \big) 
+ \omega^2 m(\boldsymbol{x})\, u.
\end{equation}
The corresponding method is called \emph{$H_k$-GenEO}, and theoretical 
estimates have been derived for the indefinite heterogeneous Helmholtz 
problem in \cite{Dolean:2024:SPW_updated}.
When Robin conditions are used, the eigenvalues of 
\eqref{DiscreteH-GenEOEigenproblemUnified} are complex but tend to 
cluster near the real axis, enabling a stable selection criterion based on
\[
\Re(\lambda) < \eta_{\text{max}}.
\]
A thorough numerical study performed in \cite{Bootland:2022:CDG} shows that $H_k$-GenEO achieves wavenumber-independent GMRES convergence and empirical frequency scaling formulae were derived. 
A robust and commonly used choice is $\eta_{\text{max}} = \frac{1}{2}$.}

\begin{remark}
Unlike the DtN coarse space, where the spectral information is extracted from the subdomain boundary and then harmonically extended into the interior, the GenEO family constructs eigenproblems that are defined \emph{throughout each subdomain}, often including the overlap region. This leads to \emph{potentially larger coarse spaces} since the spectral basis is built over the full subdomain volume. 

Although a full theoretical analysis of the GenEO methodology in the Helmholtz context is still incomplete, rigorous convergence estimates are available for certain variants. In particular, for $\Delta$-GenEO, such bounds were established in~\cite{Bootland:2022:OSM} and later refined in~\cite{Dolean:2024:ITE}. These results were recently extended to the $H_k$-GenEO construction~\cite{Dolean:2024:SPW_updated}. However, these theoretical results tend to be conservative, often overestimating the coarse space size required in practice. Numerical experiments indicate that both GenEO variants exhibit significantly better performance than what the theory predicts.
\end{remark}

\subsection{Helmholtz-harmonic coarse spaces}

We now describe the class of \emph{Helmholtz-harmonic coarse spaces}, introduced independently in~\cite{Ma:2024:TLR,Hu:2024:ANC} and later 
in a fully algebraic setting in~\cite{Nataf2024}. These spaces are built from discrete solutions of the \emph{heterogeneous} 
Helmholtz equation on local subdomains and then assembled into a global conforming approximation space using partition-of-unity weights.

\subsubsection{Harmonic GenEO coarse space}

The underlying philosophy is similar in spirit to the \emph{GenEO family}:  a local spectral decomposition is used to select relevant low-energy modes that define the coarse space.  
However, there is a key distinction: while GenEO formulates a spectral problem using the full operator (e.g.\ Laplacian or Helmholtz) and searches for dominant directions in the whole local space, \emph{harmonic coarse spaces restrict the search space to local Helmholtz-harmonic functions}. {The spectral problem is then posed over this smaller subspace typically associated with the  \emph{positive} heterogeneous Helmholtz operator \eqref{eq:posHelm}.  
}
This restriction makes the underlying spectral theory easier and the coarse space selection process \emph{real-valued}, while still encoding the Helmholtz physics through the choice of constrained subspace.

To introduce the coarse spaces described in~\cite{Ma:2024:TLR,Hu:2024:ANC}, we begin by defining the \emph{space of local Helmholtz solutions} on a subdomain $\Omega_s$:
\begin{equation}
W(\Omega_s) \assign 
\left\{
u \in V(\Omega_s):\;
a_{\Omega_s}(u, v) = 0,
\ \forall v \in V_0(\Omega_s)
\right\},
\end{equation}
where $V(\Omega_s)$ and $V_0(\Omega_s)$ are defined in \eqref{eq:vomega}.
The sesquilinear form $a_{\Omega_s}(\cdot, \cdot)$ corresponds to the Helmholtz
operator defined in~\eqref{WeakFormTermsUnified}.
Thus, $W(\Omega_s)$ consists of 
heterogeneous Helmholtz solutions with the appropriate boundary conditions.  

{To define a spectral selection on this space, we introduce the 
\emph{positive heterogeneous Helmholtz sesquilinear form}:
\begin{align}
\label{PositiveHelmholtzUnified}
b_{\Omega_s}(u,v) 
= \int_{\Omega_s} \left( a(\boldsymbol{x}) \nabla u \cdot \nabla \bar{v} 
+ \omega^2 m(\boldsymbol{x})\, u \bar{v} \right) 
\, \mathrm{d}\boldsymbol{x}.
\end{align}}
We then solve the following eigenproblem (see \cite[eq.~(3.6)]{Hu:2024:ANC} and 
\cite[eq.~(2.20)]{Ma:2024:TLR}):
\begin{align}
\label{HarmonicVariationalEigenproblemUnified}
\text{find } u \in W(\Omega_s) \quad \text{such that} \quad
b_{\Omega_s}(u,v) 
= \lambda\, b_{\Omega_s}(\Xi_s(u), \Xi_s(v)) 
\quad \forall v \in W(\Omega_s),
\end{align}
where $\Xi_s$ is the partition-of-unity operator associated with the 
subdomain $\Omega_s$.  
This problem is \emph{self-adjoint and coercive}, 
and all eigenvalues $\lambda$ are real and positive.

To handle the Helmholtz-harmonicity constraint in practice, the eigenvalue problem is recast
as a saddle point problem~\cite[eq.~(2.27)]{Ma:2024:TLR}:
find \(u \in V(\Omega_{s})\), 
\(p \in V_{0}(\Omega_{s}),\ p=0\ \text{on}\ \partial\Omega_{s}\setminus\partial\Omega\) such that

\begin{alignat}{3}\label{HarmonicVariationalEigenproblemSaddlePoint}
		b_{\Omega_s}(u,v) + {\overline{a_{\Omega_{s}}(v,p)}}
		&= \lambda b_{\Omega_s}(\Xi_s(u),\Xi_s(v)),
		&& \hspace*{10ex} \forall v \in V(\Omega_s),\\
		a_{\Omega_{s}}(u,q)
		&= 0,
		&& \hspace*{10ex}\forall q \in V_{0}(\Omega_{s}), \\
		p
		&=0,
		&& \hspace*{10ex} \text{on}\ \partial\Omega_{s}\setminus\partial\Omega.
\end{alignat} 
In matrix notations, the eigenproblem system is:
find \(\mathbf{u}\) and \(\mathbf{p}\) such that
\begin{align}\label{HarmonicEigenproblemSaddlePointMatrix}
	\begin{bmatrix}
		\widetilde{B}_{s} & {A}_{s,0}^{*}\\
		{A}_{s,0} & 0
	\end{bmatrix}
	\begin{bmatrix}
		\mathbf{u} \\ \mathbf{p}
	\end{bmatrix}
	= \lambda 
	\begin{bmatrix}
		D_{s} B_{s} D_{s} & 0\\
		0 & 0
	\end{bmatrix}
	\begin{bmatrix}
		\mathbf{u} \\ \mathbf{p}
	\end{bmatrix},
\end{align} 
where:
\begin{itemize}
  \item $A_s = R_s A R_s^T$ is the local Dirichlet matrix associated with the heterogeneous Helmholtz sesquilinear form $a_{\Omega_s}$,
  \item $A_{s,0}$ is the restriction of $A_s$ to rows not associated with boundary nodes,
  \item $B_s = R_s B R_s^T$ is the local Dirichlet matrix associated with the positive Helmholtz bilinear form $b_{\Omega_s}$,
  \item $\widetilde{B}_s$ is its Neumann counterpart,
  \item $D_s$ is the partition of unity diagonal weight matrix.
\end{itemize}
For simplicity,
the system above was written for homogeneous Dirichlet boundary conditions
\(u_{\Gamma_{D}}=0\).

After solving the eigenproblem, the \emph{coarse space basis} is formed by selecting all eigenvectors $\mathbf{u}$ associated with eigenvalues $\lambda < \tau$, where $\tau > 0$ is a user-defined threshold. These are then extended to the global space via:
\[
Z = \left\{ R_s^T D_s \mathbf{u} \right\}.
\]
Each basis function is thus globally supported (via the partition of unity), locally Helmholtz-harmonic (by construction), and selected through a stable positive definite spectral problem.

\subsubsection{Extended Harmonic GenEO coarse space}

Building on the construction of Helmholtz-harmonic coarse spaces described above, we now present a fully algebraic variant introduced in~\cite{Nataf2024}. This approach retains the same core philosophy---constructing coarse functions from locally Helmholtz-harmonic components---but introduces an \emph{extended overlapping decomposition}. 
Contrary to the previous variant, the discrete Helmholtz-harmonic constraint is
not enforced on the eigenvectors, but rather in the construction of the coarse
space matrix $Z$ thanks to a suitable operator (also present in the eigenproblem).

Specifically, for each original subdomain $\Omega_s$, an enlarged subdomain $\check{\Omega}_s$ is defined by adding at least one layer of mesh elements along its interfaces with neighbouring subdomains. This auxiliary decomposition $\{\check{\Omega}_s\}_{s=1}^N$ is used solely for the construction of the coarse space and does \emph{not} affect the first-level (fine-scale) preconditioner.

All quantities associated with the extended subdomains are denoted with a superscript $\check{\cdot}$. For example, $\check{R}_s$ is the restriction matrix to $\check{\Omega}_s$, and $\check{A}_s = \check{R}_s A \check{R}_s^T$ is the local Dirichlet matrix on $\check{\Omega}_s$.

The local eigenproblem reads: find $\check{\mathbf{u}}$ such that
\begin{align}\label{ExtendedHarmonicEigenproblem}
\widetilde{\check{B}}_{s}
\check{\mathbf{u}}
= \lambda 
\left(R_s \check{R}_s^T - \widehat{A}_s^{-1} R_s \check{R}_s^T \check{A}_s \right)^{\dagger}
D_s B_s D_s
\left(R_s \check{R}_s^T - \widehat{A}_s^{-1} R_s \check{R}_s^T \check{A}_s \right)
\check{\mathbf{u}},
\end{align}
where:
\begin{itemize}
  \item $\check{A}_s = \check{R}_s A \check{R}_s^T$ is the local Dirichlet matrix associated with the heterogeneous Helmholtz sesquilinear form $a_{\check{\Omega}_s}$ on the enlarged subdomain $\check{\Omega}_{s}$,
  \item $B_s = R_s B R_s^T$ is the local Dirichlet matrix associated with the positive Helmholtz bilinear form $b_{{\Omega}_s}$,
  \item $\widetilde{\check{B}}_s$ is the local Neumann matrix associated with the positive Helmholtz bilinear form $b_{\check{\Omega}_s}$,
  \item $\widehat{A}_s^{-1}$ corresponds to the local Robin problem on $\Omega_s$.
\end{itemize}

This eigenproblem is written under the assumption of no Dirichlet boundary conditions ($\Gamma_D = \emptyset$) for simplicity, but can accommodate them as long as the action of $\widehat{A}_s^{-1} R_s \check{R}_s^T \check{A}_s$ yields a solution satisfying the problem boundary conditions.

By construction, this eigenproblem is \emph{self-adjoint and coercive}, 
so all eigenvalues $\lambda$ are real and positive.
Moreover, while $B_s$ and $\widetilde{\check{B}}_{s}$ are typically associated
with the positive Helmholtz matrix, other SPD matrices could also be used in
their place.
We point out that the matrix on the right-hand side does not have to be
assembled in practice to solve the eigenproblems; it is enough to know the
matrix-vector product.

The coarse space basis is formed by collecting all vectors of the form
\[
Z = \left\{ R_s^T D_s 
\left( R_s \check{R}_s^T - \widehat{A}_s^{-1} R_s \check{R}_s^T \check{A}_s \right)
\check{\mathbf{u}} \;:\; \lambda < \tau \right\}
\]
for some user-defined threshold $\tau > 0$.
By construction, each resulting vector lies in the global finite element space
and corresponds to a \emph{local discrete Helmholtz solution} on $\Omega_s$
thanks to the presence of the operator
$R_s \check{R}_s^T - \widehat{A}_s^{-1} R_s \check{R}_s^T \check{A}_s$,
and multiplied by the partition of unity $D_{s}$.

\subsection{Summary and comparison of spectral coarse spaces}

Having described the main families of spectral coarse spaces---namely the DtN, GenEO-type, and Helmholtz-harmonic constructions---we now highlight their connections and comparative features.
A summary of this comparison is given in Table~\ref{tab:coarse_space_comparison}.

\paragraph{Conceptual links}

Despite differing technical constructions, all
approaches share a common spectral philosophy: they extract low-energy
components via localised eigenproblems, then assemble them into globally
supported basis functions using a partition of unity.

The \textbf{DtN} coarse space is the most geometrically minimal, using
interface-local eigenfunctions extended via Helmholtz harmonic continuation.
In this sense, it belongs to the broader class of \textbf{Helmholtz-harmonic
coarse spaces}.
It may also be seen as a precursor to the \textbf{GenEO} family, which replaces
boundary-local eigenproblems with volumetric ones defined over subdomains or
overlapping regions.
The \textbf{harmonic GenEO} spaces, in turn, reinterpret the GenEO framework by
restricting the usual self-adjoint positive definite eigenproblems to the
subspace of local Helmholtz solutions.
The constraint to solve in the Helmholtz-harmonic subspace can be imposed using
an augmented eigenproblem.
Finally, the \textbf{extended-harmonic} coarse spaces involve self-adjoint
positive definite eigenproblems similarly to the GenEO family, but posed on an
extended subdomain and with a suitable operator weighting the matrix in the
eigenproblem.
The resulting eigenvectors are projected in the Helmholtz-harmonic subspace of
the (non-extended) subdomain.

{
\paragraph{Complexity comparison}
The main costs of two-level domain decomposition preconditioners are:
(i) construction of the coarse space (local eigenproblems and global assembly),
(ii) factorisation of the coarse operator, and
(iii) repeated local and coarse solves during GMRES iterations.

For discretisations with a fixed number of points per wavelength (\(hk =
\mathcal{O}(1)\)) as often happens in engineering practice, the number of unknowns per
subdomain scales like \(\mathcal{O}(k^d)\) in \(d\) dimensions, while the number of
unknowns on the subdomain boundary scales like \(\mathcal{O}(k^{d-1})\).
As a result the cost of constructing the coarse space and its dimension differ
between the different methods:
\begin{itemize}
	\item \textbf{DtN:} the eigenproblem is posed on the interface (although the
		operator features a local volume inverse of size \(\mathcal{O}(k^d)\))
		and as a result, the number of basis functions is bounded by the number
		of interface unknowns in \(\mathcal{O}(k^{d-1})\). 
	\item \textbf{GenEO-type:} the eigenproblem is volumetric hence leads to
		\(\mathcal{O}(k^d)\) candidate functions per subdomain.
		To ensure wavenumber-independent convergence, the coarse operator
		(although well below this limit in practice) can become large, leading
		to an important factorisation cost.
	\item \textbf{Harmonic:} The saddle-point eigenproblem is twice the
		number of unknowns in the subdomain, hence of size
		\(\mathcal{O}(k^d)\), but the retained modes need to solve the
		Helmholtz problem, hence their number is necessarily bounded by the number
		of boundary unknowns in \(\mathcal{O}(k^{d-1})\).
	\item \textbf{Extended harmonic:} The eigenproblem is posed in the
		extended volume, still in \(\mathcal{O}(k^d)\), but the elements in the
		coarse space are Helmholtz solutions by construction (up to the
		partition of unity), hence their number is also necessarily bounded by
		the number of boundary unknowns in \(\mathcal{O}(k^{d-1})\).
\end{itemize}
At runtime, each GMRES iteration requires local subdomain solves (local
matrices are factored in the setup phase) typically via sparse direct solvers
and a global coarse correction (\(\mathcal{O}(N_c^2)\) after factorisation). In
principle, these scalings suggest coarse space dimensions increase dramatically
with the wavenumber, but in practice our numerical results show milder growth,
with iteration counts nearly independent of \(k\) once a suitable coarse space
is included.}

\begin{table}[ht]
\footnotesize
\centering
\caption{Comparison of spectral coarse space families for the Helmholtz problem}
\label{tab:coarse_space_comparison}
\renewcommand{\arraystretch}{1.2}
\begin{tabular}{|l|c|c|c|c|}
\hline
\textbf{Property} & \textbf{DtN} & \textbf{GenEO-type} & \textbf{Harmonic} & {\textbf{Extended harmonic}} \\
\hline
\textbf{Eigenfunctions} & Interface $\Gamma_s$ & Subdomain $\Omega_s$ (or overlap) & Subdomain $\Omega_s$ & {Extended subdomain $\check{\Omega}_s$} \\
\hline
\textbf{Spectral domain} & Interface & Full volume & Helmholtz-harmonic subspace & Full extended volume \\
\hline
\textbf{Eigenvalues} & complex & complex & real positive & real positive \\
\hline
\textbf{Maximum dimension} & $\mathcal{O}(|\Gamma_s|)$ & $\mathcal{O}(|\Omega_s|)$ & $\mathcal{O}(|\Gamma_s|)$ & $\mathcal{O}(|\Gamma_s|)$ \\
\hline
\textbf{Theory} & None (empirically) & Available for some variants & Available & Available\\
\hline
\textbf{Adapted to high $k$} & Yes & Yes (\(H_{k}\)-GenEO) & Yes & Yes \\
\hline
\textbf{Notable references} & \cite{Conen:2014:ACS} & \cite{Spillane:2014:ARC}, \cite{Bootland:2022:OSM}, \cite{Dolean:2024:ITE}, \cite{Dolean:2024:SPW_updated} & \cite{Hu:2024:ANC}, \cite{Ma:2024:TLR} & \cite{Nataf2024} \\
\hline
\end{tabular}
\end{table}

\paragraph{Rationale for numerical comparison} Given the subtle theoretical distinctions and the absence of definitive complexity bounds in the indefinite Helmholtz regime, it is difficult to draw conclusions purely from the analytical form of the coarse spaces. The actual performance depends on practical implementation details, the choice of eigenvalue threshold, and the frequency regime. Therefore, a numerical comparison is essential to:
\begin{itemize}
    \item assess the robustness of each coarse space across a range of frequencies and heterogeneous media;
    \item quantify the trade-off between coarse space dimension and solver convergence;
    \item validate whether the observed empirical behaviour aligns with or diverges from existing theoretical predictions.
\end{itemize}

In the next section, we provide such a numerical study, focusing on both homogeneous and heterogeneous Helmholtz problems with increasing complexity.

\section{Numerical assessment and comparison}
\label{sec:ComparativeNumericalStudies}

The purpose of this section is to evaluate the performance and robustness of the spectral coarse spaces introduced in Section~4 in realistic simulation settings. While these methods differ in their mathematical formulation and spectral construction, it is often difficult to draw firm conclusions about their practical efficiency based solely on theoretical arguments. In particular, the coarse space dimension, conditioning effects, and preconditioner performance depend intricately on the frequency, domain geometry, heterogeneity, and the choice of eigenvalue threshold.

To overcome this, we conduct a series of numerical experiments designed to probe the behaviour of each coarse space under increasing complexity
evaluated within a common two-level ORAS preconditioning framework, using GMRES as the iterative solver. In all cases, the GMRES solver {is not restarted} 
and the tolerance is \(10^{-6}\) on the residual, with a right-preconditioning. The overlap is minimal (unless stated otherwise) with a symmetry constraint with respect to the
interface, implying two layers of cells in the overlap. Each test case is analysed through a combination of strong and weak scaling experiments, studying the influence of:
\begin{itemize}
  \item subdomain diameter (in wavelength),
  \item coarse space size (global and per subdomain),
  \item eigenvalue threshold $\tau$,
  \item and overlap and partition of unity choices.
\end{itemize}

Our goal is to assess both efficiency (number of GMRES iterations) and scalability (coarse space dimension and distribution), and to offer a fair comparison between the different coarse space strategies under practically relevant conditions.

{\paragraph{Implementation details} The FreeFEM implementation of the benchmarks is done within the ffddm framework, a set of parallel FreeFEM scripts implementing Schwarz domain decomposition methods. The ffddm documentation is available on the FreeFem.org web page (see~\cite{FFD:Tournier:2019}). The ffddm framework relies on Message Passing Interface (MPI) parallelism. As is usually done in domain decomposition methods, we assign one subdomain per MPI process.

The first step is to decompose the computational domain into overlapping subdomains. Unless otherwise stated, the automatic graph partitioner Metis~\cite{Karypis:98:METIS} is used, which produces a non-overlapping decomposition of the set of mesh elements while minimising interfaces between subdomains and conserving good load-balancing. The overlapping decomposition is then obtained by adding successive layers of elements to reach the desired size of overlap. The setup of the one-level preconditioner then consists in assembling and factorising the local matrices in each subdomain in parallel. This is done using the sparse direct solver MUMPS~\cite{amestoy2001fully}. The two-level preconditioner for each method is then built by first solving the corresponding local eigenvalue problems in each subdomain with SLEPc~\cite{hernandez2005ssf}, and finally assembling and factorising the coarse space operator in a distributed manner by MUMPS using a few cores (ranging from 1 to 144 depending on the size of the coarse space).

During the solution phase, each application of the preconditioner involves solving linear systems with local subdomain matrices (first level) and with the coarse space matrix (second level), which is done by forward-backward substitution using the factorisations computed during the setup phase.

A FreeFEM script comparing all methods for the homogeneous Helmholtz equation in a square is available at \url{https://github.com/FreeFem/FreeFem-sources/blob/develop/examples/ffddm/Helmholtz-2d-all.edp}.}

\subsection{Numerical simulations in a square: homogeneous and heterogeneous problems}
The numerical tests in this section are structured as follows:
\begin{enumerate}
  \item \textbf{Homogeneous test case.} A baseline benchmark on a unit square domain with Robin boundary conditions and a single interior point source. This allows us to measure the scaling of coarse space dimension and iteration count with respect to the number of subdomains and frequency, in a basic controlled setting.
  
  \item \textbf{Heterogeneous test case.} A layered medium with strong contrasts in wave speed is introduced to study robustness under medium heterogeneity and its impact on coarse space size and convergence.
\end{enumerate}

\clearpage
\subsubsection{Homogeneous problem}
\begin{figure}[H]
	\centering
	\begin{tikzpicture}[scale=1.75]
		% Grid
		\draw[step=0.5cm,gray!50,very thin, shift={(-1,-1)}] (0,0) grid (2,2);
		\draw[gray!50,very thin] (-1,-1) -- (1,1);
		\draw[gray!50,very thin] (-1,1) -- (1,-1);
		\draw[gray!50,very thin] (-1,0) -- (0,1) -- (1,0) -- (0,-1) -- cycle;
		% Boundaries
		\fill[gray,opacity=0.2] (-1,1) -- (-1,-1) -- (1,-1) -- (1,1) -- cycle;
		\draw (-1,1) -- (-1,-1) -- (1,-1) -- (1,1) -- cycle;
		% Boundary labels
		\draw (1.05,0) node[above,rotate=-90] {Robin};
		\draw (-1.05,0) node[above,rotate=90] {Robin};
		\draw (0,-1.05) node[below] {Robin};
		\draw (0,1.05) node[above] {Robin};
		% Coordinate labels
		\draw (-1.05,1) node[left] {$y = 1$}
		(-1.05,-1) node[left] {$y = 0$}
		(1,-1.05) node[below] {$x = 1$}
		(-1,-1.05) node[below] {$x = 0$};
		% Point source
		\draw[fill=black] (0,0) circle (0.025) node[below] {source $f$};
	\end{tikzpicture}
	\hfill
	{\includegraphics[width=0.33\textwidth]{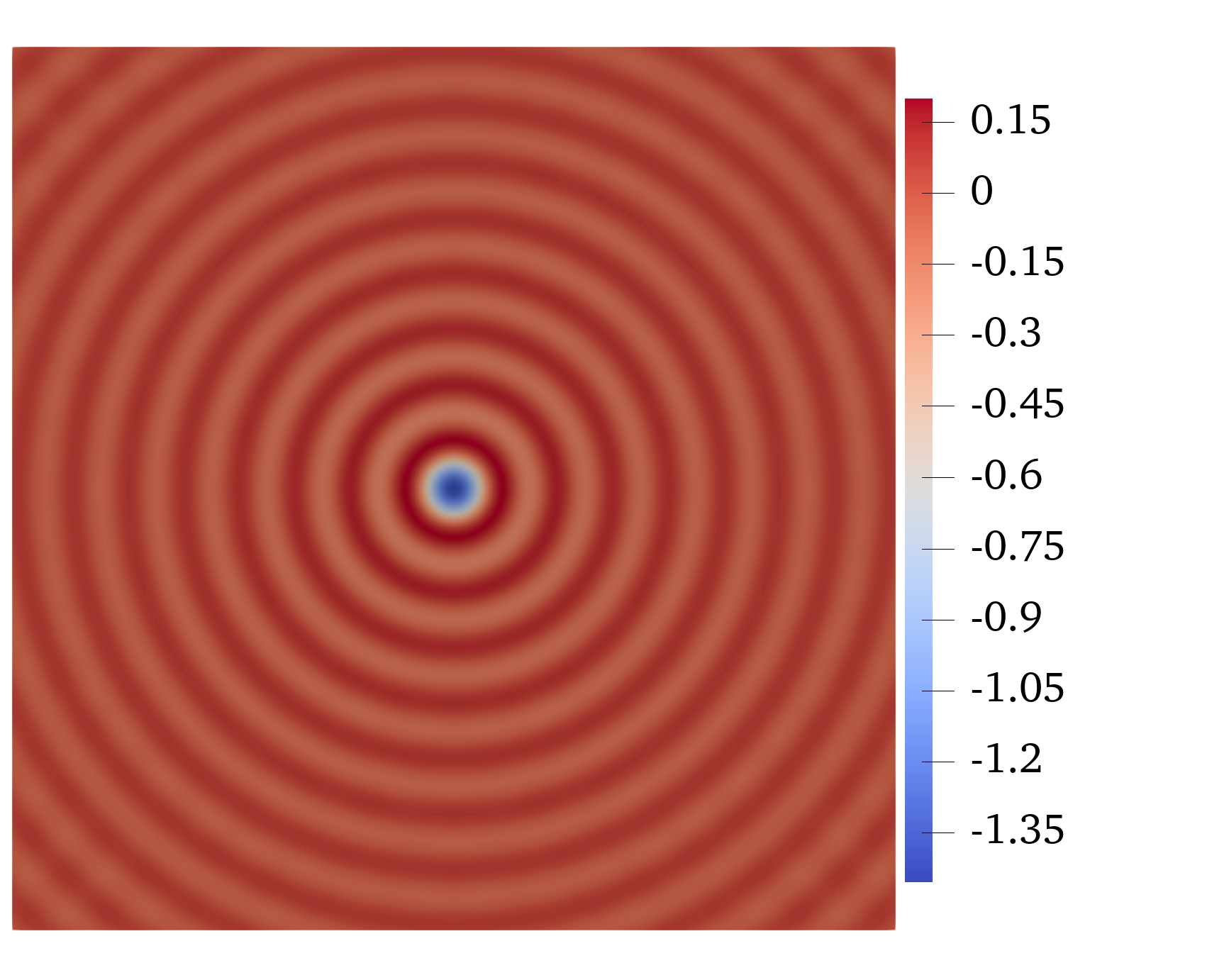}}
	{\includegraphics[width=0.33\textwidth]{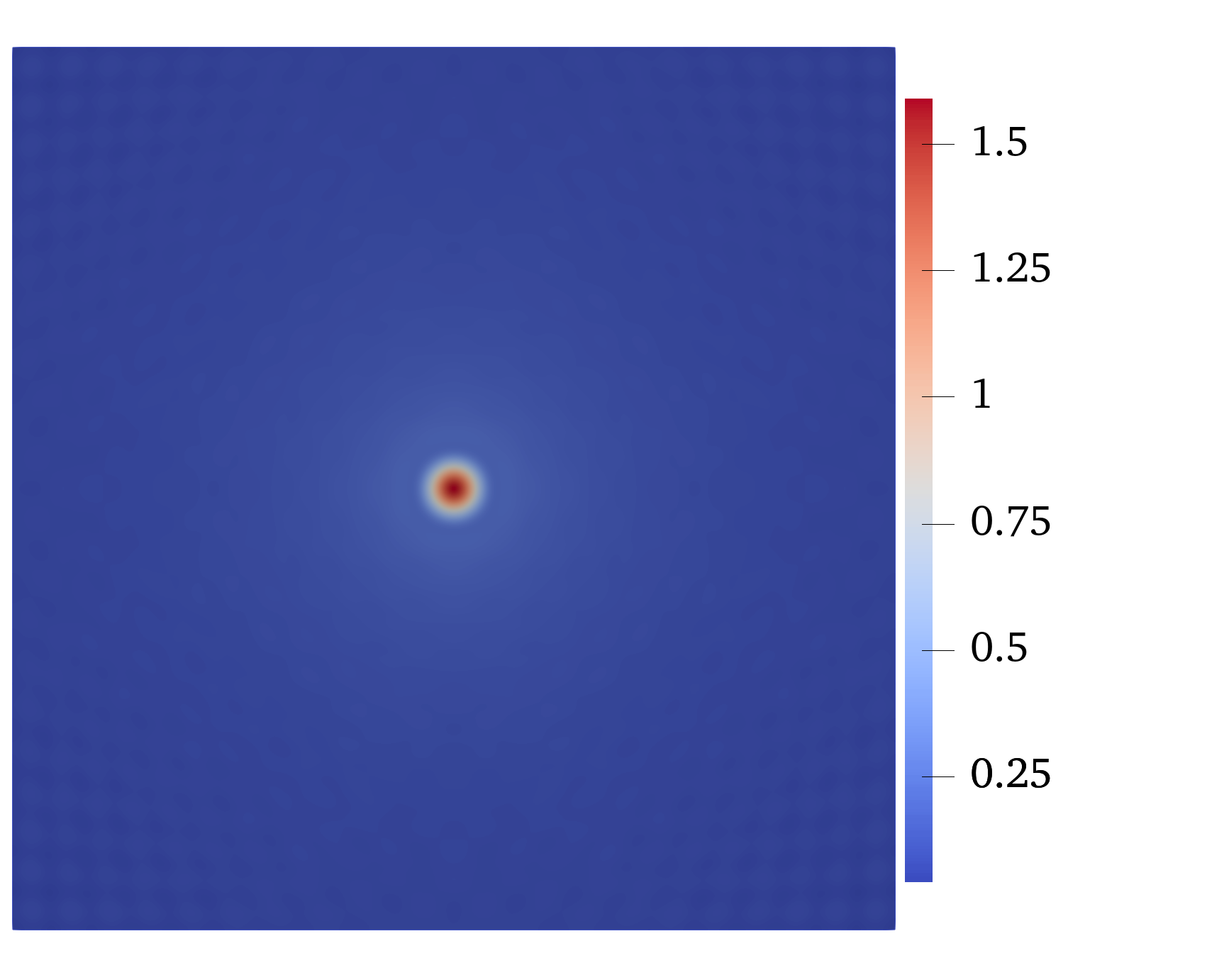}}
	\caption{
		Schematic of the 2D wave guide model problem with example triangular mesh.
 		Real part (left) and modulus (right) of the total field
 		for the homogeneous media test case with \(\omega=100\).
    }\label{Fig:WaveGuide2D}
\end{figure}
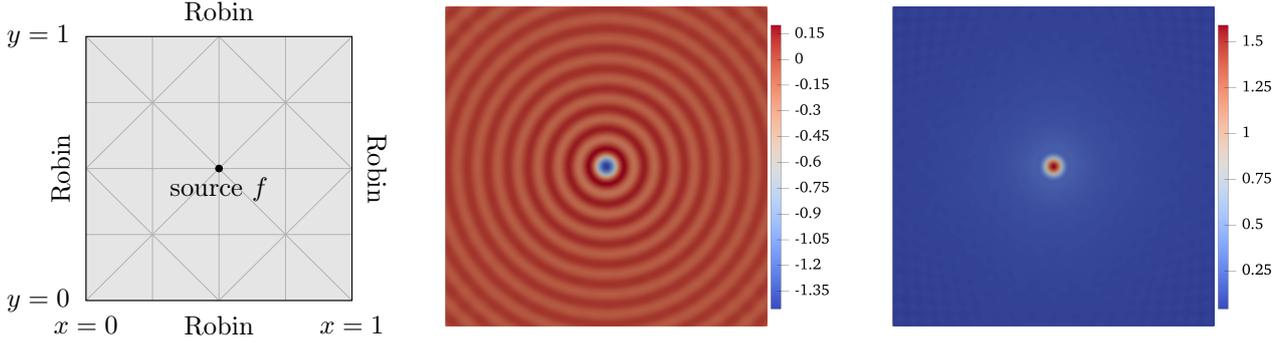

In this part we base the numerical experiments on the model problem in 2D, defined on the unit square $\Omega = (0,1)^2$. We impose Robin boundary conditions on all sides of the domain. A point source is located in the centre of the domain at $(\frac{1}{2},\frac{1}{2})$ and provides the forcing function $f$.
The point source is numerically modelled by a Gaussian function:
\(f(x,y)=10^4 \exp(-10^3[(x-\frac{1}{2})^2+(y-\frac{1}{2})^2])\).
A schematic of this model problem is found in Figure~\ref{Fig:WaveGuide2D}. 

To discretise the problem, we triangulate $\Omega$ using a Cartesian grid with spacing $h$ and alternating diagonals to form a simplicial mesh.
{We consider constant coefficients $a(\boldsymbol{x}) \equiv 1$ and} 
$m(\boldsymbol{x}) \equiv 1$. 
The local wavenumber \( k = \omega \sqrt{\frac{m}{a}} {= \omega} \) is then constant
and the wavelength used to measure geometrical parameters is \(\lambda_{\min}=2\pi/k\).
The discrete problem \eqref{LinearSystemUnified} is assembled using a $\mathbb{P}_2$ finite element approximation on this mesh.  
To mitigate the \emph{pollution effect}, we choose the angular frequency  $\omega$ (or equivalently the wavenumber $k$) and the mesh size $h$ simultaneously so that the dimensionless quantity $k h$ remains sufficiently small.  
In practice, this is enforced by fixing a minimum number of grid points per wavelength $\lambda_{\min}$; here, we ensure that 
$k h \lesssim 0.5$, which corresponds to at least 10 points per wavelength $\lambda_{\min}$.

\begin{figure}[H]
    \centering
	{\includegraphics[width=0.49\textwidth]{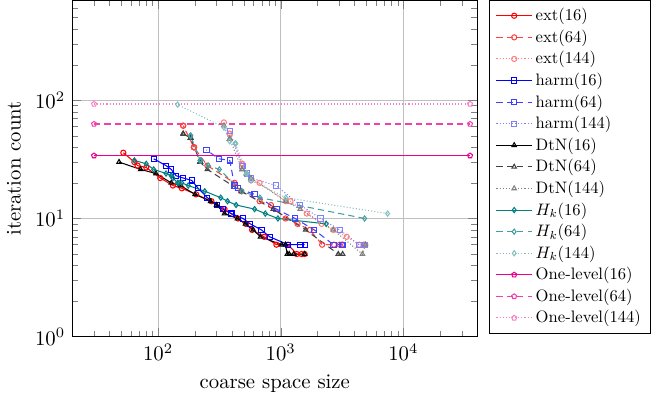}}
	{\includegraphics[width=0.49\textwidth]{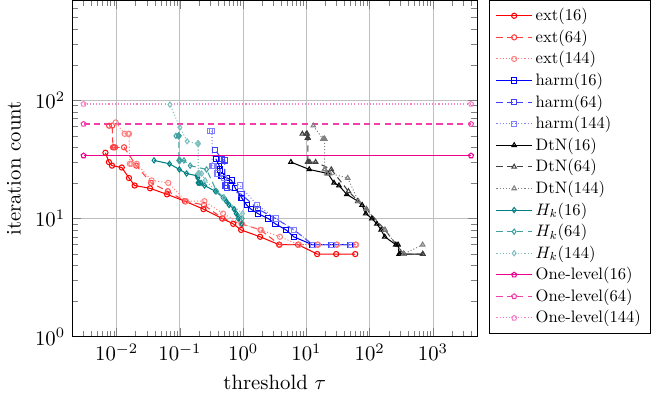}}
	{\includegraphics[width=0.49\textwidth]{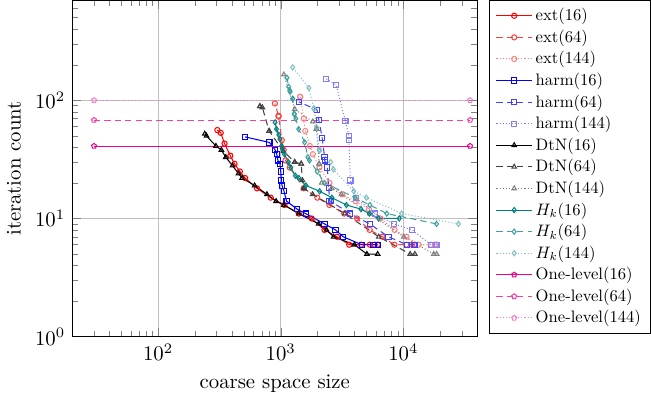}}
	{\includegraphics[width=0.49\textwidth]{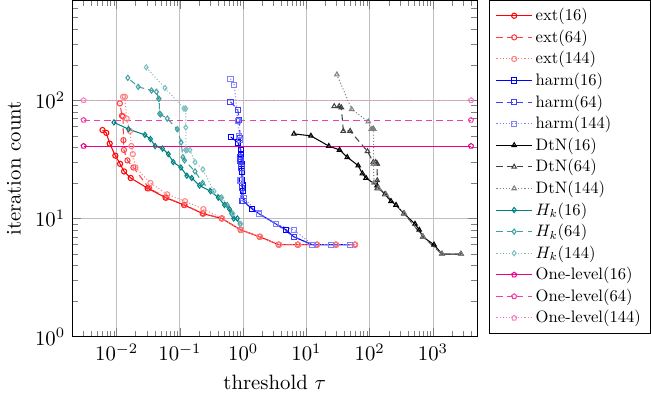}}
    \caption{
		Influence of the coarse space size (left) and threshold choice (right)
		on the iteration count in strong scaling for the
		homogeneous media test case with \(\omega=20\) (top) and \(\omega=100\)
		(bottom).
        {The number in brackets indicates the number of subdomains.}
    }\label{fig:strong_scaling_vscssize_homog}
\end{figure}

\begin{table}[H]
	\scriptsize
	\begin{center}
\begin{tabular}{ cccccc|c|ccc|ccc|ccc|ccc }
\multicolumn{6}{c|}{} & \multicolumn{1}{c|}{1lvl} & \multicolumn{3}{c|}{ext} & \multicolumn{3}{c|}{harm} & \multicolumn{3}{c|}{DtN} & \multicolumn{3}{c}{$H_k$} \\
    $L$ [$\lambda_{\min}$] & $N$ & $n$ & $H$ [$\lambda_{\min}$] & $n_{s}$ & $n_{s}^{\partial\Omega_{s}}$ & It & It & CS & CS$_{s}$ & It & CS & CS$_{s}$ & It & CS & CS$_{s}$ & It & CS & CS$_{s}$ \\
\hline
 4.5 & 16 & 231361 & 1.1 & 15373 & 490 & 34 & 5 & 1352 & 84 & 6 & 1148 & 72 & 5 & 1132 & 71 & 9 & 2352 & 147 \\
 4.5 & 36 & 231361 & 0.8 & 7108 & 331 & 47 & 6 & 1528 & 42 & 6 & 1904 & 53 & 5 & 2102 & 58 & 10 & 3576 & 99 \\
 4.5 & 64 & 231361 & 0.6 & 4156 & 251 & 63 & 6 & 2176 & 34 & 6 & 2690 & 42 & 5 & 2946 & 46 & 10 & 4832 & 76 \\
 4.5 & 100 & 231361 & 0.5 & 2762 & 203 & 77 & 6 & 2784 & 28 & 6 & 3504 & 35 & 5 & 3802 & 38 & 11 & 6094 & 61 \\
 4.5 & 144 & 231361 & 0.4 & 1990 & 171 & 93 & 6 & 4368 & 30 & 6 & 4366 & 30 & 5 & 4636 & 32 & 11 & 7440 & 52 \\
\hline
 9.0 & 16 & 519841 & 2.3 & 33853 & 730 & 36 & 6 & 1360 & 85 & 6 & 1724 & 108 & 5 & 1906 & 119 & 9 & 3504 & 219 \\
 9.0 & 36 & 519841 & 1.5 & 15455 & 491 & 48 & 6 & 2298 & 64 & 6 & 2864 & 80 & 5 & 3116 & 87 & 10 & 5304 & 147 \\
 9.0 & 64 & 519841 & 1.1 & 8926 & 371 & 63 & 6 & 3200 & 50 & 6 & 3996 & 62 & 5 & 4342 & 68 & 10 & 7136 & 112 \\
 9.0 & 100 & 519841 & 0.9 & 5863 & 299 & 73 & 6 & 5140 & 51 & 6 & 5140 & 51 & 5 & 5572 & 56 & 10 & 8962 & 90 \\
 9.0 & 144 & 519841 & 0.8 & 4177 & 251 & 92 & 6 & 6480 & 45 & 6 & 6478 & 45 & 5 & 6834 & 47 & 10 & 10896 & 76 \\
\hline
 13.5 & 16 & 1442401 & 3.4 & 92413 & 1210 & 39 & 5 & 3346 & 209 & 6 & 2856 & 178 & 5 & 3156 & 197 & 10 & 3316 & 207 \\
 13.5 & 36 & 1442401 & 2.3 & 41748 & 811 & 52 & 6 & 4748 & 132 & 7 & 4768 & 132 & 6 & 4182 & 116 & 10 & 8760 & 243 \\
 13.5 & 64 & 1442401 & 1.7 & 23866 & 611 & 67 & 6 & 5316 & 83 & 6 & 6684 & 104 & 5 & 7164 & 112 & 10 & 11744 & 184 \\
 13.5 & 100 & 1442401 & 1.4 & 15520 & 491 & 82 & 6 & 6874 & 69 & 6 & 8576 & 86 & 5 & 9148 & 91 & 10 & 14760 & 148 \\
 13.5 & 144 & 1442401 & 1.1 & 10950 & 411 & 97 & 6 & 8448 & 59 & 6 & 10460 & 73 & 5 & 11186 & 78 & 10 & 17808 & 124 \\
\hline
 18.0 & 16 & 2076481 & 4.5 & 132493 & 1450 & 39 & 6 & 3404 & 213 & 6 & 4288 & 268 & 7 & 3782 & 236 & 10 & 6960 & 435 \\
 18.0 & 36 & 2076481 & 3.0 & 59695 & 971 & 51 & 6 & 5676 & 158 & 7 & 5728 & 159 & 6 & 5020 & 139 & 10 & 10488 & 291 \\
 18.0 & 64 & 2076481 & 2.3 & 34036 & 731 & 67 & 6 & 8000 & 125 & 7 & 8028 & 125 & 6 & 7044 & 110 & 9 & 14048 & 220 \\
 18.0 & 100 & 2076481 & 1.8 & 22077 & 587 & 81 & 6 & 8180 & 82 & 6 & 10344 & 103 & 5 & 10950 & 110 & 10 & 17596 & 176 \\
 18.0 & 144 & 2076481 & 1.5 & 15537 & 491 & 98 & 6 & 10074 & 70 & 6 & 13992 & 97 & 6 & 13340 & 93 & 10 & 21264 & 148 \\
\hline
 22.5 & 16 & 3690241 & 5.6 & 234253 & 1930 & 41 & 6 & 3616 & 226 & 6 & 4570 & 286 & 5 & 5042 & 315 & 10 & 6280 & 392 \\
 22.5 & 36 & 3690241 & 3.8 & 105188 & 1291 & 54 & 6 & 6036 & 168 & 6 & 7612 & 211 & 5 & 8198 & 228 & 9 & 13944 & 387 \\
 22.5 & 64 & 3690241 & 2.8 & 59776 & 971 & 68 & 6 & 8452 & 132 & 6 & 10692 & 167 & 5 & 11358 & 177 & 9 & 18656 & 292 \\
 22.5 & 100 & 3690241 & 2.3 & 38647 & 779 & 84 & 6 & 14940 & 149 & 7 & 13778 & 138 & 6 & 14526 & 145 & 9 & 23362 & 234 \\
 22.5 & 144 & 3690241 & 1.9 & 27110 & 651 & 100 & 6 & 13344 & 93 & 6 & 16796 & 117 & 5 & 17694 & 123 & 9 & 28176 & 196 \\
\end{tabular}
\end{center}

    \caption{
		Strong scaling experiment for the homogeneous media test case.
    }\label{tab:strong_scaling_homog}
\end{table}

\begin{figure}[H]
    \centering
	{\includegraphics[width=0.48\textwidth]{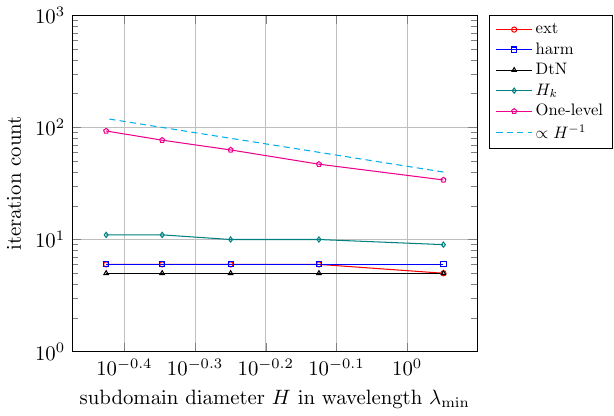}}
	{\includegraphics[width=0.48\textwidth]{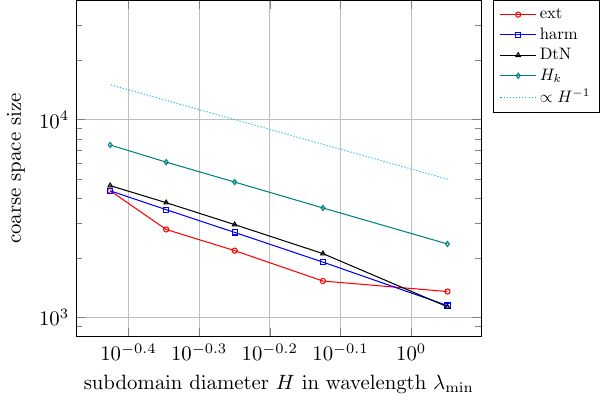}}
	{\includegraphics[width=0.48\textwidth]{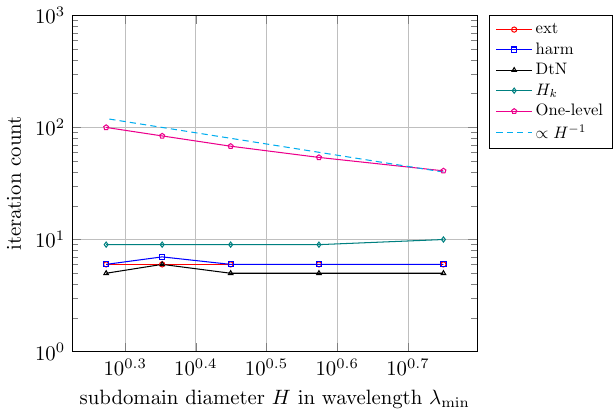}}
	{\includegraphics[width=0.48\textwidth]{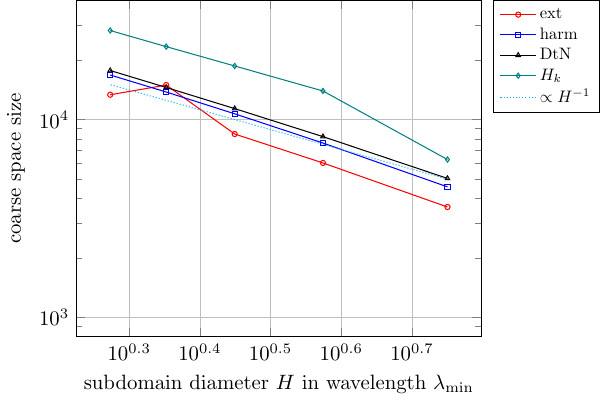}}
    \caption{
		Influence of the subdomain diameter on the iteration count (left) and
		coarse space size (right) in strong scaling for the homogeneous media
		test case \(\rho=10\) with \(\omega=20\) (top) and \(\omega=100\)
		(bottom).
    }\label{fig:strong_scaling_vsdiameter_homog}
\end{figure}

\clearpage
\subsubsection{Heterogeneous problem}

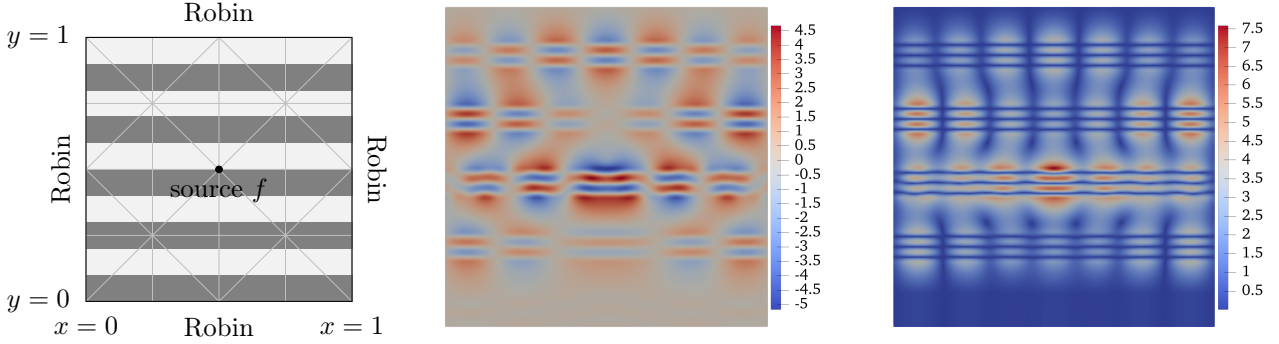
\begin{figure}[H]
	\centering
	\begin{tikzpicture}[scale=1.75]
		%grey bands
		\fill[gray,opacity=0.1] (-1,1) -- (1,1) -- (1,0.8) -- (-1,0.8) -- cycle;
		\fill[gray,opacity=1.0] (-1,0.6) -- (1,0.6) -- (1,0.8) -- (-1,0.8) -- cycle;
		\fill[gray,opacity=0.1] (-1,0.6) -- (1,0.6) -- (1,0.4) -- (-1,0.4) -- cycle;
		\fill[gray,opacity=1.0] (-1,0.2) -- (1,0.2) -- (1,0.4) -- (-1,0.4) -- cycle;
		\fill[gray,opacity=0.1] (-1,0.2) -- (1,0.2) -- (1,0.0) -- (-1,0.0) -- cycle;
		\fill[gray,opacity=1.0] (-1,-0.2) -- (1,-0.2) -- (1,0.0) -- (-1,0.0) -- cycle;
		\fill[gray,opacity=0.1] (-1,-0.2) -- (1,-0.2) -- (1,-0.4) -- (-1,-0.4) -- cycle;
		\fill[gray,opacity=1.0] (-1,-0.6) -- (1,-0.6) -- (1,-0.4) -- (-1,-0.4) -- cycle;
		\fill[gray,opacity=0.1] (-1,-0.6) -- (1,-0.6) -- (1,-0.8) -- (-1,-0.8) -- cycle;
		\fill[gray,opacity=1.0] (-1,-1) -- (1,-1) -- (1,-0.8) -- (-1,-0.8) -- cycle;
		% Grid
		\draw[step=0.5cm,gray!50,very thin, shift={(-1,-1)}] (0,0) grid (2,2);
		\draw[gray!50,very thin] (-1,-1) -- (1,1);
		\draw[gray!50,very thin] (-1,1) -- (1,-1);
		\draw[gray!50,very thin] (-1,0) -- (0,1) -- (1,0) -- (0,-1) -- cycle;
		% Boundaries
		\draw (-1,1) -- (-1,-1) -- (1,-1) -- (1,1) -- cycle;
		% Boundary labels
		\draw (1.05,0) node[above,rotate=-90] {Robin};
		\draw (-1.05,0) node[above,rotate=90] {Robin};
		\draw (0,-1.05) node[below] {Robin};
		\draw (0,1.05) node[above] {Robin};
		% Coordinate labels
		\draw (-1.05,1) node[left] {$y = 1$}
			  (-1.05,-1) node[left] {$y = 0$}
			  (1,-1.05) node[below] {$x = 1$}
			  (-1,-1.05) node[below] {$x = 0$};
		% Point source
		\draw[fill=black] (0,0) circle (0.025) node[below] {source $f$};
	\end{tikzpicture}
	\hfill
	{\includegraphics[width=0.33\textwidth]{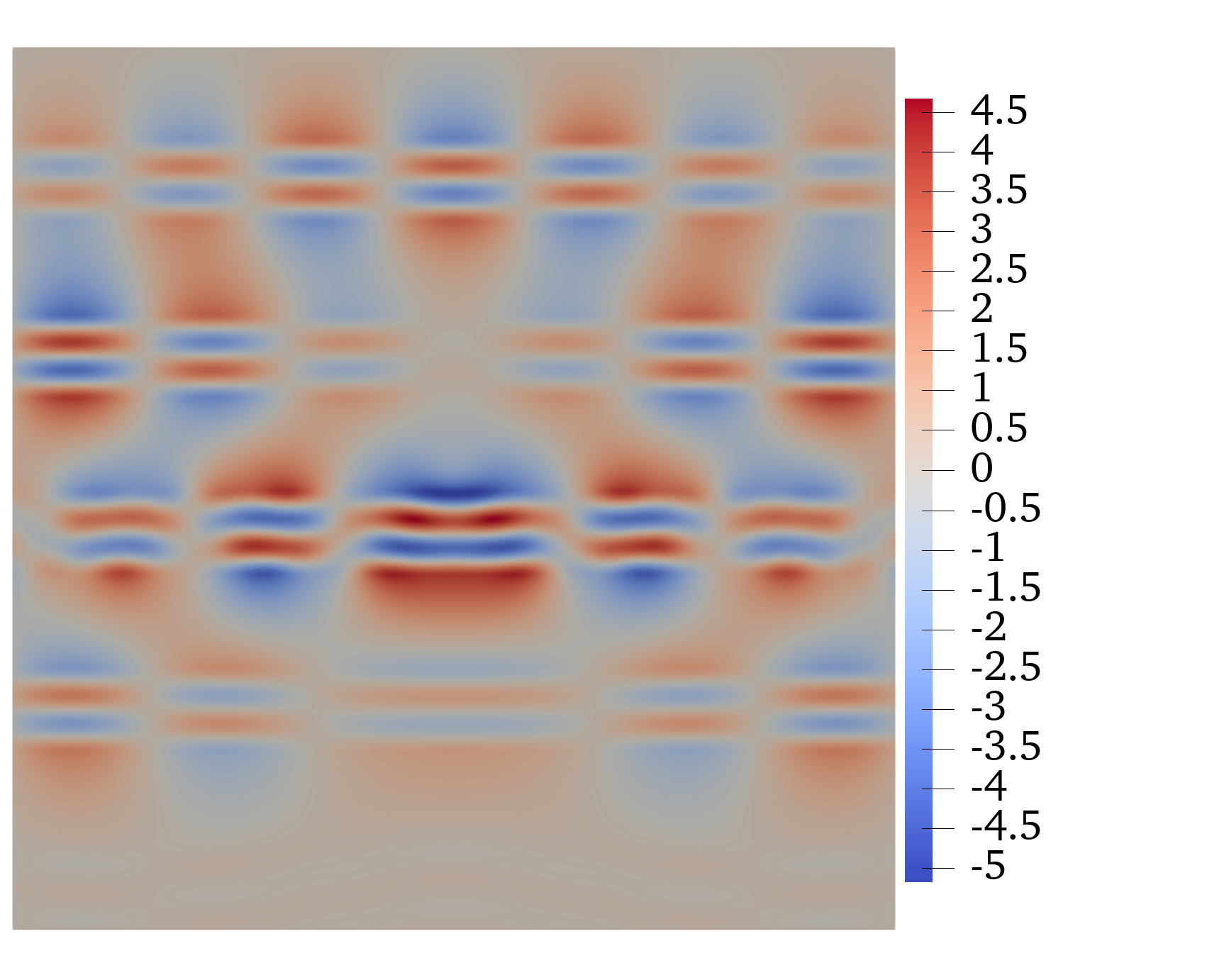}}
	{\includegraphics[width=0.33\textwidth]{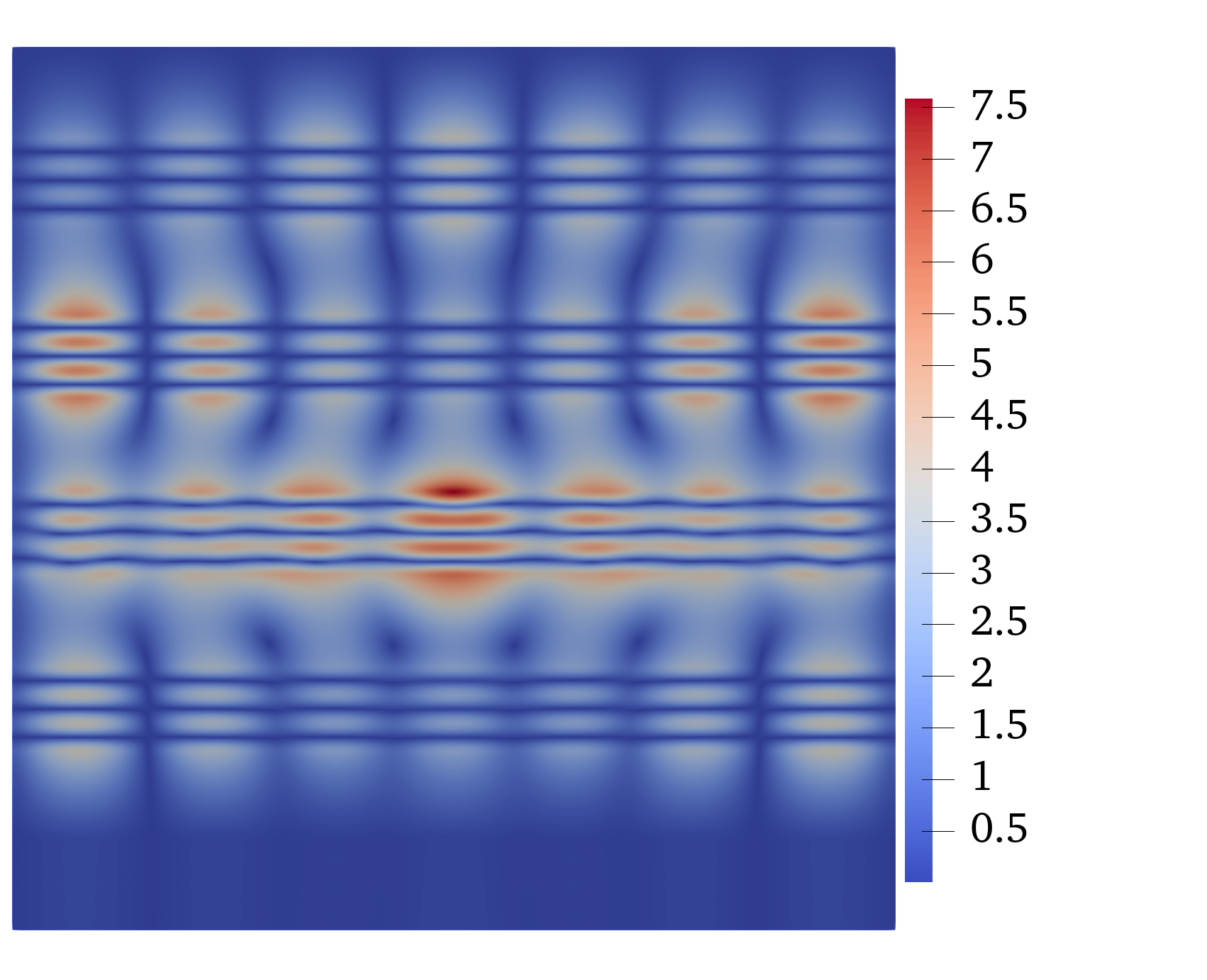}}
    \caption{
		Piecewise constant layer profiles for the wave speed
		$c(\boldsymbol{x})$ (left).
		For the darkest shade, $c(\boldsymbol{x}) = 1$, while for the lightest
		shade, $c(\boldsymbol{x}) = \rho$, with $\rho$ being the contrast
		factor.
		Real part (middle) and modulus (right) of the total field for the
		heterogeneous media test case: \(\rho=10\) and \(\omega=100\).
    }\label{fig:solutions_heterogeneous}
\end{figure}

The next numerical experiments aim to study the robustness of the method with
respect to medium heterogeneity.
For this purpose, we consider a layered medium in the unit square.
More precisely, we use an \emph{alternating layer} configuration (see
Figure~\ref{fig:solutions_heterogeneous}) in which the heterogeneity is
introduced  through the material parameter $m(\boldsymbol{x})$.
In the present setting, we keep $a(\boldsymbol{x}) \equiv 1$ everywhere and
the mass coefficient is defined as
$m(\boldsymbol{x}) = c(\boldsymbol{x})^{-2}$,  
where $c(\boldsymbol{x})$ is the spatially varying wave speed.  
The medium alternates between two constant wave speeds:
\(
c(\boldsymbol{x}) \in \{ 1, \ \rho \},
\)
where $\rho > 1$ is a contrast parameter controlling the strength of the heterogeneity.  The local wavenumber is then given by
\(
k(\boldsymbol{x}) 
= \omega \sqrt{\frac{m(\boldsymbol{x})}{a(\boldsymbol{x})}} 
= \frac{\omega}{c(\boldsymbol{x})},
\)
with $\omega > 0$ the angular frequency.  
In our numerical experiments, we vary $\omega$ while keeping $\rho$ fixed  unless otherwise stated.
Geometric parameters are measured using the minimal wavelength \(\lambda_{\min} = 2\pi/\omega\) corresponding to \(c(\boldsymbol{x})=1\).
An example of a numerical solution for the heterogeneous layered medium 
is shown in Figure~\ref{fig:solutions_heterogeneous}.

\begin{figure}[H]
    \centering
	{\includegraphics[width=0.49\textwidth]{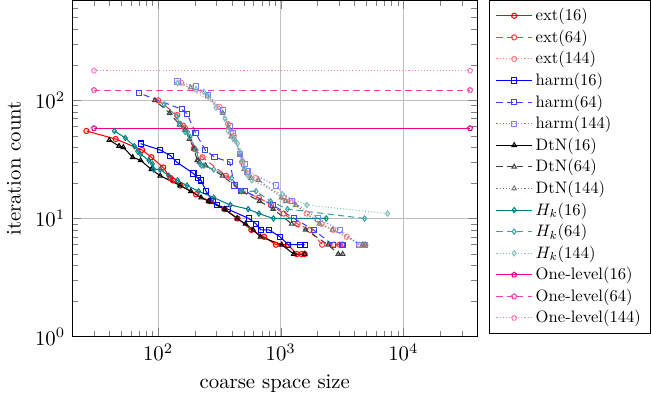}}
	{\includegraphics[width=0.49\textwidth]{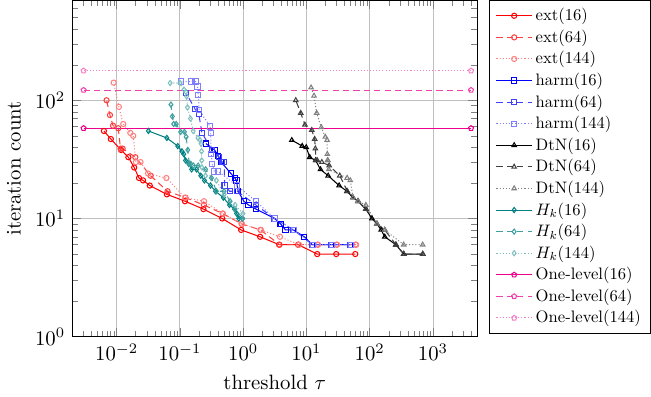}}
	{\includegraphics[width=0.49\textwidth]{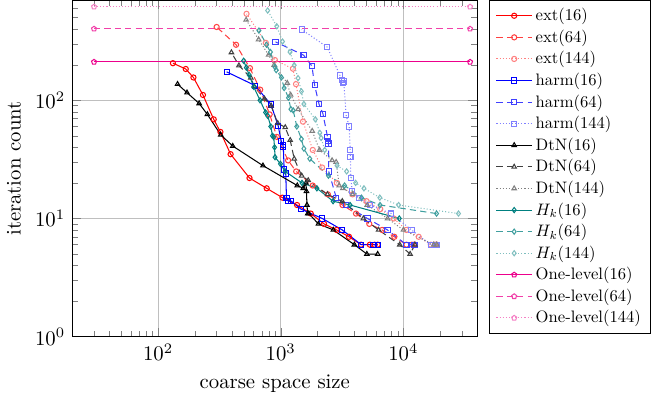}}
	{\includegraphics[width=0.49\textwidth]{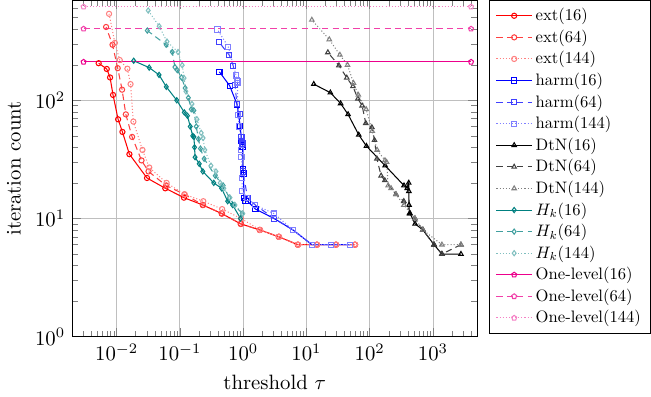}}
    \caption{
		Influence of the coarse space size (left) and threshold choice (right)
		on the iteration count in strong scaling for the
		heterogeneous media test case \(\rho=10\) with \(\omega=20\) (top) and
		\(\omega=100\) (bottom).
    }\label{fig:strong_scaling_vscssize_layers_cmax10}
\end{figure}

\begin{table}[H]
	\scriptsize
	\begin{center}
\begin{tabular}{ cccccc|c|ccc|ccc|ccc|ccc }
\multicolumn{6}{c|}{} & \multicolumn{1}{c|}{1lvl} & \multicolumn{3}{c|}{ext} & \multicolumn{3}{c|}{harm} & \multicolumn{3}{c|}{DtN} & \multicolumn{3}{c}{$H_k$} \\
    $L$ [$\lambda_{\min}$] & $N$ & $n$ & $H$ [$\lambda_{\min}$] & $n_{s}$ & $n_{s}^{\partial\Omega_{s}}$ & It & It & CS & CS$_{s}$ & It & CS & CS$_{s}$ & It & CS & CS$_{s}$ & It & CS & CS$_{s}$ \\
\hline
 4.5 & 16 & 231361 & 1.1 & 15373 & 490 & 58 & 5 & 1352 & 84 & 6 & 1148 & 72 & 5 & 1276 & 80 & 10 & 872 & 54 \\
 4.5 & 36 & 231361 & 0.8 & 7108 & 331 & 94 & 6 & 1528 & 42 & 6 & 1904 & 53 & 5 & 2100 & 58 & 10 & 3576 & 99 \\
 4.5 & 64 & 231361 & 0.6 & 4156 & 251 & 122 & 6 & 2176 & 34 & 6 & 2652 & 41 & 5 & 2946 & 46 & 10 & 4832 & 76 \\
 4.5 & 100 & 231361 & 0.5 & 2762 & 203 & 154 & 6 & 3500 & 35 & 6 & 3504 & 35 & 5 & 3802 & 38 & 11 & 6094 & 61 \\
 4.5 & 144 & 231361 & 0.4 & 1990 & 171 & 178 & 6 & 4368 & 30 & 6 & 4356 & 30 & 6 & 4636 & 32 & 11 & 7440 & 52 \\
\hline
 9.0 & 16 & 519841 & 2.3 & 33853 & 730 & 119 & 6 & 1710 & 107 & 6 & 1712 & 107 & 5 & 1906 & 119 & 10 & 3504 & 219 \\
 9.0 & 36 & 519841 & 1.5 & 15455 & 491 & 178 & 6 & 2860 & 79 & 6 & 2864 & 80 & 6 & 2538 & 70 & 11 & 5304 & 147 \\
 9.0 & 64 & 519841 & 1.1 & 8926 & 371 & 231 & 6 & 4032 & 63 & 6 & 4644 & 73 & 6 & 3572 & 56 & 11 & 7136 & 112 \\
 9.0 & 100 & 519841 & 0.9 & 5863 & 299 & 294 & 6 & 5140 & 51 & 6 & 5924 & 59 & 6 & 5572 & 56 & 11 & 8962 & 90 \\
 9.0 & 144 & 519841 & 0.8 & 4177 & 251 & 327 & 6 & 6409 & 45 & 6 & 7137 & 50 & 6 & 5764 & 40 & 12 & 10896 & 76 \\
\hline
 13.5 & 16 & 1442401 & 3.4 & 92413 & 1210 & 84 & 5 & 3345 & 209 & 6 & 2838 & 177 & 5 & 3154 & 197 & 10 & 5808 & 363 \\
 13.5 & 36 & 1442401 & 2.3 & 41748 & 811 & 117 & 5 & 5372 & 149 & 6 & 4732 & 131 & 6 & 4177 & 116 & 10 & 8760 & 243 \\
 13.5 & 64 & 1442401 & 1.7 & 23866 & 611 & 157 & 7 & 6657 & 104 & 8 & 6660 & 104 & 13 & 3956 & 62 & 10 & 11744 & 184 \\
 13.5 & 100 & 1442401 & 1.4 & 15520 & 491 & 205 & 6 & 6848 & 68 & 6 & 8560 & 86 & 5 & 9148 & 91 & 10 & 14760 & 148 \\
 13.5 & 144 & 1442401 & 1.1 & 10950 & 411 & 198 & 6 & 8388 & 58 & 6 & 10460 & 73 & 5 & 11186 & 78 & 10 & 17808 & 124 \\
\hline
 18.0 & 16 & 2076481 & 4.5 & 132493 & 1450 & 120 & 6 & 3395 & 212 & 6 & 3404 & 213 & 8 & 3778 & 236 & 10 & 6960 & 435 \\
 18.0 & 36 & 2076481 & 3.0 & 59695 & 971 & 160 & 6 & 5676 & 158 & 7 & 5676 & 158 & 9 & 6165 & 171 & 10 & 10488 & 291 \\
 18.0 & 64 & 2076481 & 2.3 & 34036 & 731 & 218 & 6 & 7964 & 124 & 7 & 7946 & 124 & 11 & 8562 & 134 & 10 & 14048 & 220 \\
 18.0 & 100 & 2076481 & 1.8 & 22077 & 587 & 294 & 6 & 10180 & 102 & 6 & 10196 & 102 & 5 & 10948 & 109 & 10 & 17596 & 176 \\
 18.0 & 144 & 2076481 & 1.5 & 15537 & 491 & 288 & 6 & 12561 & 87 & 6 & 13992 & 97 & 6 & 13340 & 93 & 10 & 21264 & 148 \\
\hline
 22.5 & 16 & 3690241 & 5.6 & 234253 & 1930 & 211 & 6 & 4520 & 282 & 6 & 4532 & 283 & 5 & 5039 & 315 & 10 & 9264 & 579 \\
 22.5 & 36 & 3690241 & 3.8 & 105188 & 1291 & 306 & 6 & 6014 & 167 & 6 & 7544 & 210 & 5 & 8198 & 228 & 10 & 13944 & 387 \\
 22.5 & 64 & 3690241 & 2.8 & 59776 & 971 & 402 & 6 & 10592 & 166 & 6 & 10592 & 166 & 5 & 11357 & 177 & 11 & 18656 & 292 \\
 22.5 & 100 & 3690241 & 2.3 & 38647 & 779 & 553 & 9 & 10826 & 108 & 8 & 15407 & 154 & 7 & 14526 & 145 & 10 & 23362 & 234 \\
 22.5 & 144 & 3690241 & 1.9 & 27110 & 651 & 618 & 6 & 16754 & 116 & 6 & 16704 & 116 & 6 & 17693 & 123 & 11 & 28176 & 196 \\
\end{tabular}
\end{center}

    \caption{
		Strong scaling experiment for the heterogeneous media test case 
		\(\rho=10\).
    }\label{tab:strong_scaling_layers_cmax10}
\end{table}

\begin{figure}[H]
    \centering
	{\includegraphics[width=0.48\textwidth]{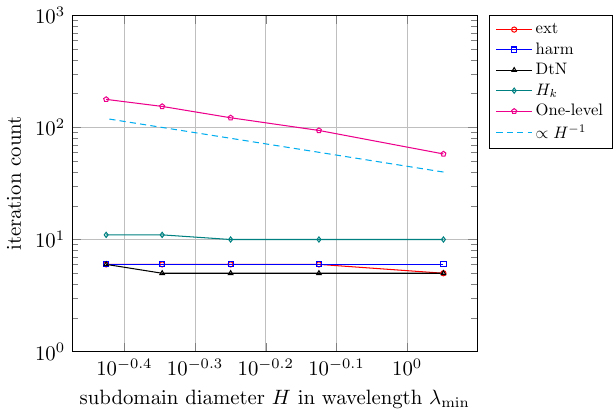}}
	{\includegraphics[width=0.48\textwidth]{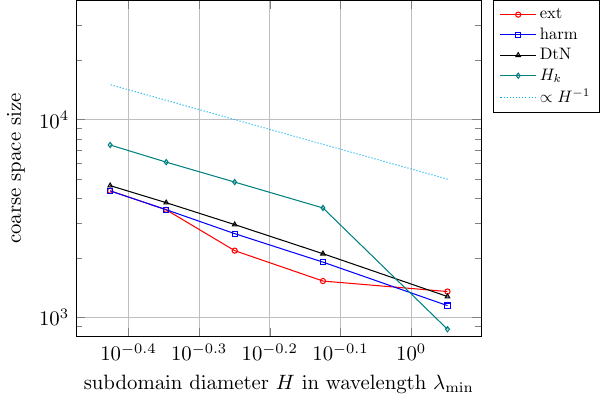}}
	{\includegraphics[width=0.48\textwidth]{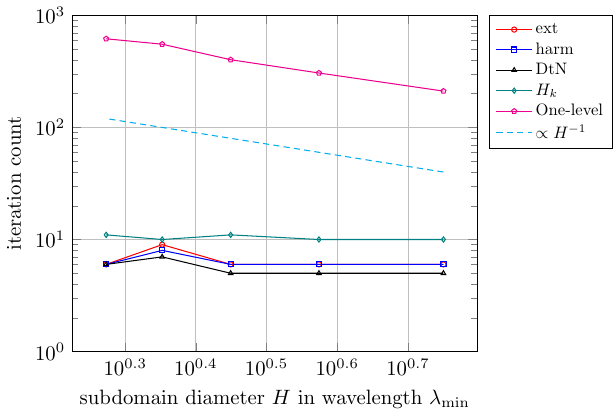}}
	{\includegraphics[width=0.48\textwidth]{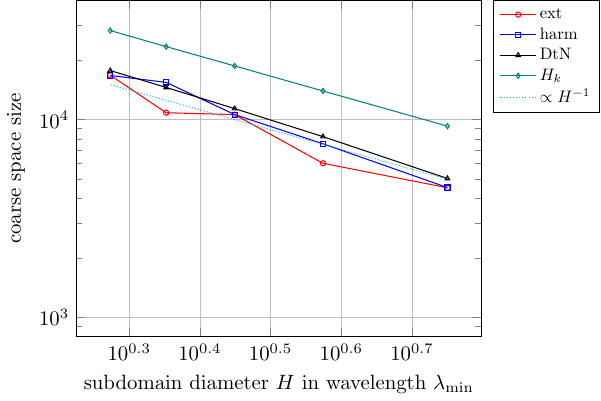}}
    \caption{
		Influence of the subdomain diameter on the iteration count (left) and
		coarse space size (right) in strong scaling for the heterogeneous media
		test case \(\rho=10\) with \(\omega=20\) (top) and \(\omega=100\)
		(bottom).
    }\label{fig:strong_scaling_vsdiameter_layers_cmax10}
\end{figure}

\subsubsection{Results}

The full results are provided in Tables \ref{tab:strong_scaling_homog} and  \ref{tab:strong_scaling_layers_cmax10}.
The entries are the diameter of the domain \(L\) measured in the smallest wavelength \(\lambda_{\min} = 2\pi\min_{\boldsymbol{x}}(c(\boldsymbol{x}))/\omega\) present in the problem, the total number of degrees of freedom \(n\), the averaged diameter of a subdomain \(H\) measured in wavelength \(\lambda_{\min}\), the averaged
number \(n_{s}\) of degrees of freedom in a subdomain, the averaged number \(n_{s}^{\partial\Omega_{s}}\) of degrees of freedom on the boundary of a subdomain, the number \(N\) of subdomains, the total size of the coarse space
(CS), the averaged number of contributions of a subdomain to the coarse space (CS$_{s}$) and the number of GMRES iterations to reach the tolerance~(It).
{The methods being} compared are the baseline one-level method (1lvl): extended
(ext) and harmonic coarse spaces (harm) on one side with the classical the
Dirichlet-to-Neumann (DtN) and \(H_{k}\)-Geneo methods ($H_k$) on the other
side.

Besides, we report, for the smallest and largest frequencies tested, the number
of GMRES iterations to reach the required tolerance as a function of the total
coarse space size and choice of threshold parameter in Figures
\ref{fig:strong_scaling_vscssize_homog} and
\ref{fig:strong_scaling_vscssize_layers_cmax10}.
The number of subdomains is provided in brackets in the legend entries.
The results reported in the table correspond to the parameters with the minimum
number of iterations obtained for each method.

The number of iterations of GMRES to reach the required tolerance as a function of the subdomain diameter (measured in wavelength) is reported
Figures \ref{fig:strong_scaling_vsdiameter_homog} and \ref{fig:strong_scaling_vsdiameter_layers_cmax10}.

\subsubsection{Assessment and comparative conclusions for the square test cases}

We now draw conclusions from the strong scaling tests performed on the 2D square homogeneous and heterogeneous test cases. 
From both Tables~\ref{tab:strong_scaling_homog}--\ref{tab:strong_scaling_layers_cmax10} and Figures~\ref{fig:strong_scaling_vscssize_homog}--\ref{fig:strong_scaling_vsdiameter_layers_cmax10}, several trends emerge consistently across frequencies and levels of heterogeneity.
General comments on one-level and two-level methods are:
\begin{itemize}
	\item The performance of the one-level method deteriorates
		with increasing frequency or with heterogeneities.
	\item In contrast, all two-level methods are scalable and very robust with
		respect to the frequency or in the presence of heterogeneity.
		As a result, they dramatically outperform the one-level baseline.
	\item These properties are however achieved only for large enough coarse
		spaces, as the iteration count of two-level methods correlates strongly
		with coarse space size.
	\item In particular, robustness is achieved with a coarse space size that
		increases with increasing number of subdomains and increasing
		frequency.
\end{itemize}
Comparing the different two-level methods, we remark that:
\begin{itemize}
	\item The smallest number of iterations is achieved with the DtN coarse
		spaces, closely followed by the harmonic and extended-harmonic coarse
		spaces, and then by the \(H_{k}\)-GenEO coarse spaces.
	\item To achieve this, the DtN coarse spaces are however slightly larger
		than the harmonic and extended-harmonic ones, while the \(H_{k}\)-GenEO
		coarse spaces are the largest.
	\item In a context with a tight coarse space size budget, the DtN and
		extended-harmonic coarse spaces are the most efficient, especially at
		large frequencies.
	\item A sensible choice for the eigenvalue threshold \(\tau\) to ensure low
		iteration count is
		in the range \(5-10\) for extended-harmonic, 
		around \(10\) for harmonic,
		in the range \(500-1000\) for DtN and
		in the range \(0.8-1\) for \(H_{k}\)-GenEO. 
\end{itemize}

\subsection{A test case from medical imaging}

\paragraph{Description of the problem}

We consider the problem of plane wave scattering by randomly positioned, 
penetrable micro-reflectors in two dimensions.  
This setting models acoustic wave interaction with soft biological tissues, 
as explored in recent quantitative ultrasound imaging 
studies~\cite{Garnier2023,Garnier2025}.
We fix a constant \emph{reference medium} with coefficients 
\(
a_0 > 0, \ m_0 > 0,
\)
and define the \emph{reference wave speed} and the corresponding \emph{reference wavenumber}
and \emph{reference wavelength}
\[
c_0 = \sqrt{\frac{a_0}{m_0}},\quad
k_0 = \omega \sqrt{\frac{m_0}{a_0}} = \frac{\omega}{c_0},\quad
\lambda_0 = \frac{2\pi}{k_0}.
\]
All geometric parameters in this test case will be expressed in units of \(\lambda_{0}\).

The incident wave is a plane wave of the form
\(
u_i(\mathbf{x}) = \exp\!\left( \imath\, k_0 \, \mathbf{d} \cdot \mathbf{x} \right),
\)
\(
\mathbf{d} \in \mathbb{R}^2, \ |\mathbf{d}| = 1,
\)
where $\mathbf{d}$ is the propagation direction.
The total field $u$ satisfies the heterogeneous Helmholtz equation
\[
- \nabla \cdot \big( a(\mathbf{x}) \nabla u \big) 
- \omega^2 m(\mathbf{x}) \, u = 0,
\quad \text{in } \mathbb{R}^2,
\]
with an outgoing radiation condition on the scattered field 
$u_s = u - u_i$.
The coefficients $a(\mathbf{x})$ and $m(\mathbf{x})$ describe the \(N_{r}\in\mathbb{N}\) penetrable 
micro-reflectors:
\begin{equation}
\label{eq:coeffs_microreflectors_k0}
a(\mathbf{x}) = 
\boldsymbol{1}_{\Omega \setminus \bigcup_{j=1}^{N_r} B_\epsilon(\mathbf{x}_j)} 
+ \sum_{j=1}^{N_r} a_j \, \boldsymbol{1}_{B_\epsilon(\mathbf{x}_j)}, 
\quad
c_0^{2} m(\mathbf{x}) = 
\boldsymbol{1}_{\Omega \setminus \bigcup_{j=1}^{N_r} B_\epsilon(\mathbf{x}_j)} 
+ \sum_{j=1}^{N_r} m_j \, \boldsymbol{1}_{B_\epsilon(\mathbf{x}_j)},
\end{equation}
where $B_\epsilon(\mathbf{x}_j)$ is a disk of radius
$\epsilon = \lambda_0 / 4$ centred at $\mathbf{x}_j$,
$a_j$ and $m_j$ are drawn from a normal distribution with mean $1$
and standard deviation $\sigma = 0.1$,
the disk centres $\mathbf{x}_j$ are generated by a Matèrn hardcore spatial
process.
In the background medium $(a_0, m_0)$ we have $k(\mathbf{x}) \equiv k_0$, 
while inside each inclusion $k(\mathbf{x})$ deviates due to the perturbed 
material parameters.

We truncate the unbounded domain using a perfectly matched layer (PML) 
following~\cite{Galkowski2024}.  
Let $\chi \in C^\infty_{\mathrm{comp}}(\mathbb{R}^2)$ be a smooth cut-off 
function supported around the heterogeneous region that is identically one on an open connected set including the support of \((1-m)(1-a)\). We introduce the compactly supported unknown
\[
\widetilde{u} \assign u - (1 - \chi) u_i,
\]
which satisfies
\[
- \nabla \cdot \big( a(\mathbf{x}) \nabla \widetilde{u} \big) 
- \omega^2 m(\mathbf{x}) \, \widetilde{u} 
= 2\, \nabla \chi \cdot \nabla u_i + u_i \, \Delta \chi,
\quad \text{in } \mathbb{R}^2.
\]
We use a PML implementation adapted from~\cite{Bermudez2004}.

\begin{figure}[H]
    \centering
    \includegraphics[width=0.32\textwidth]{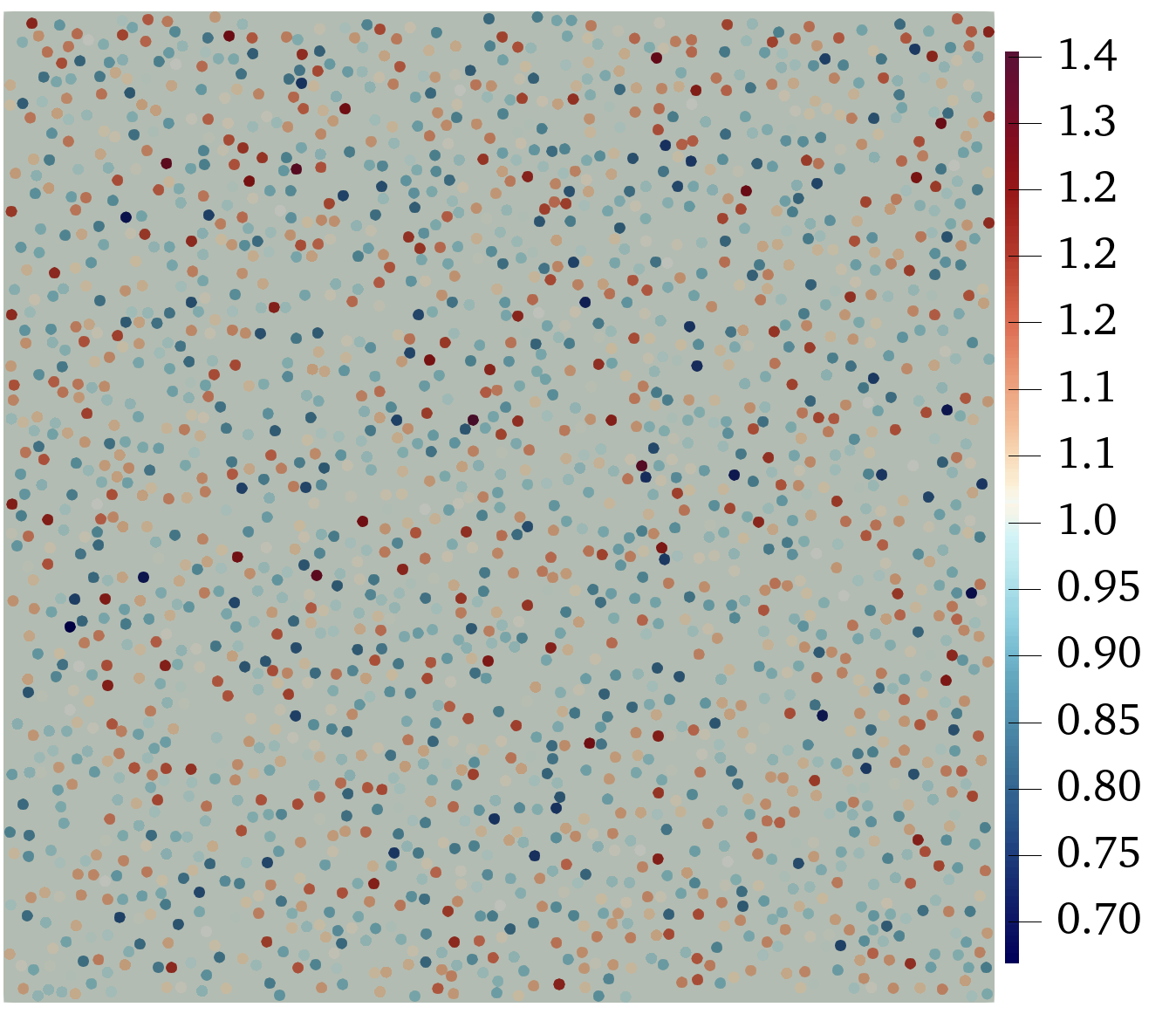}
    \includegraphics[width=0.32\textwidth]{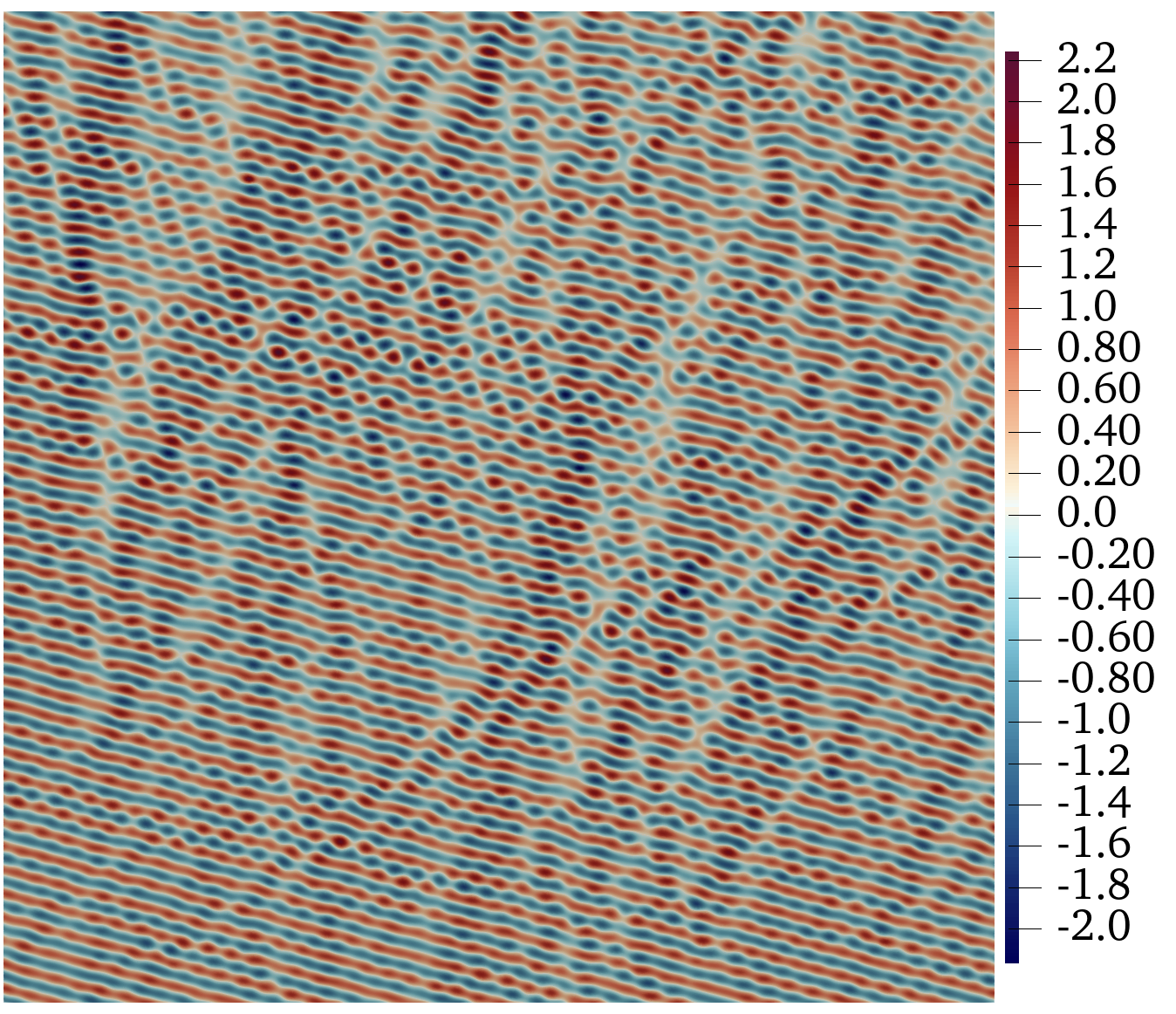}
    \includegraphics[width=0.32\textwidth]{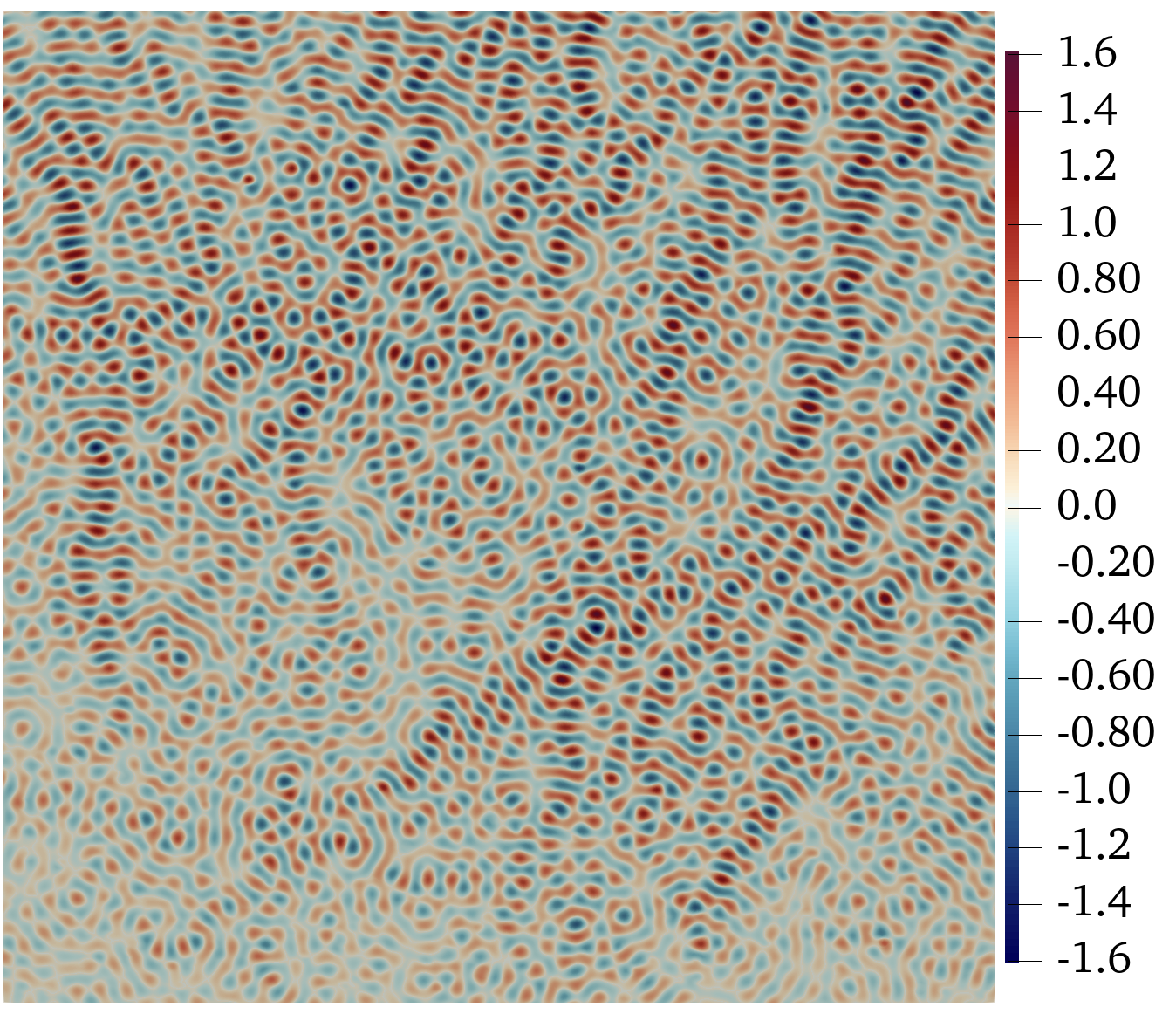}
    \caption{Values of the coefficient \(c_{0}^{2}m\) (left), real part of total field \(\Re(u)\) (middle), and real part of scattered field \(\Re(u_s)\) (right).}
    \label{fig:solutions_imaging}
\end{figure}

A snapshot of the medium and computed fields is shown in 
Figure~\ref{fig:solutions_imaging} for 
$N_r = 2040$ reflectors in a domain of size 
$47.8 \lambda_0 \times 47.8 \lambda_0$.
The presence of the micro-reflectors imposes a strong constraint on the mesh
size corresponding to about \(22\) points per wavelength \(\lambda_{0}\).
We use \(\mathbb{P}_{2}\) Lagrange finite elements.

\paragraph{Strong scaling test}

\begin{figure}[H]
    \centering
	{\includegraphics[width=0.49\textwidth]{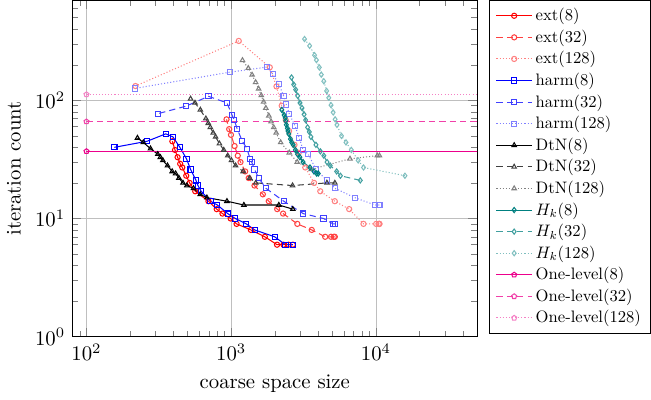}}
	{\includegraphics[width=0.49\textwidth]{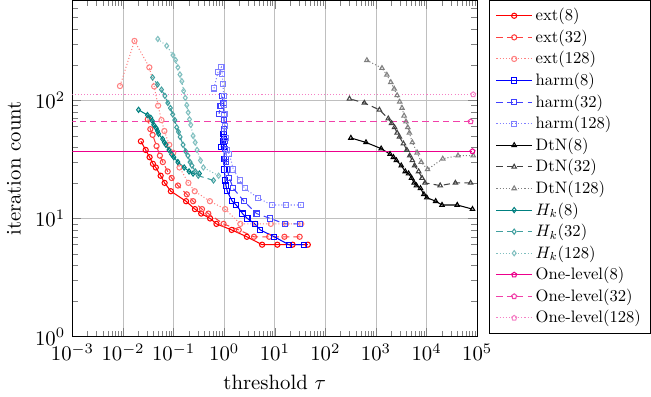}}
    \caption{
		Influence of the coarse space size (left) and threshold choice (right)
		on the iteration count in strong scaling for the imaging
		test case.
    }\label{fig:strong_scaling_vscssize_imaging}
\end{figure}

\begin{figure}[H]
    \centering
	{\includegraphics[width=0.49\textwidth]{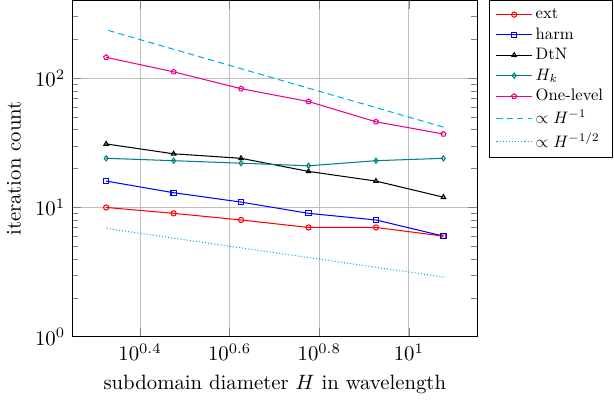}}
	{\includegraphics[width=0.49\textwidth]{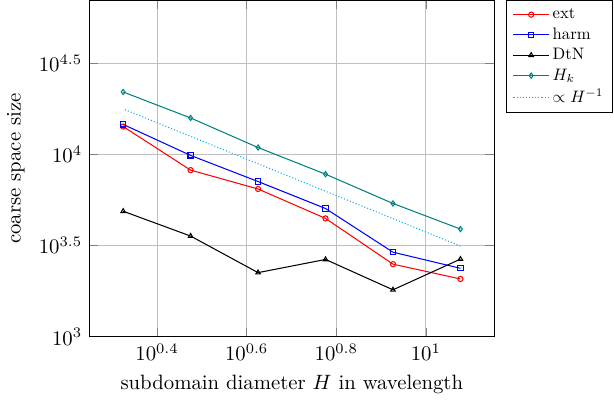}}
    \caption{
        Influence of the subdomain diameter on the iteration count in strong scaling 
		for the imaging test case.
    }\label{fig:strong_scaling_vsdiameter_imaging}
\end{figure}

We first provide results for a strong scaling test, increasing the number of
subdomains for a fixed problem size.
The full domain has size \(23.9\lambda_{0} \times 23.9\lambda_{0}\) with \(1213761\) DOFs.
The full results are provided in Table~\ref{tab:strong_scaling_imaging}.
The entries are the same as the ones previously described, except we use
the reference wavelength \(\lambda_{0}\) to measure lengths.
Besides, we report the number of iterations of GMRES to reach the required
tolerance as a function of the total coarse space size and choice of threshold
parameter in
Figure~\ref{fig:strong_scaling_vscssize_imaging}.
The number of subdomains is provided in brackets in the legend entries.
The results reported in the table correspond to the parameters with the minimum
number of iterations obtained for each method, which explains the difference
for the DtN approach.
The number of iterations of GMRES to reach the required tolerance
as a function of the subdomain diameter (measured in wavelength) is reported
Figure~\ref{fig:strong_scaling_vsdiameter_imaging}.

\paragraph{Weak scaling test}

\begin{figure}[H]
    \centering
	{\includegraphics[width=0.49\textwidth]{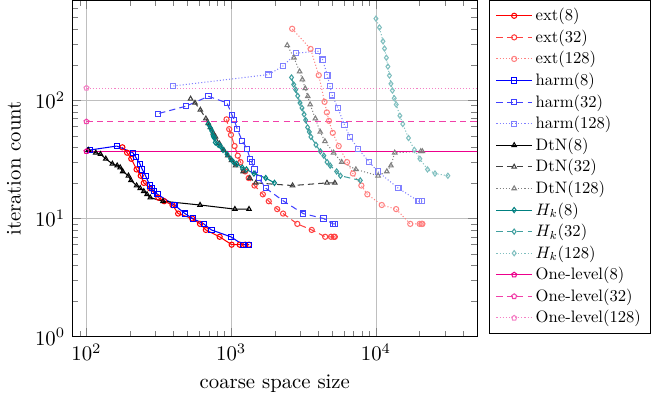}}
	{\includegraphics[width=0.49\textwidth]{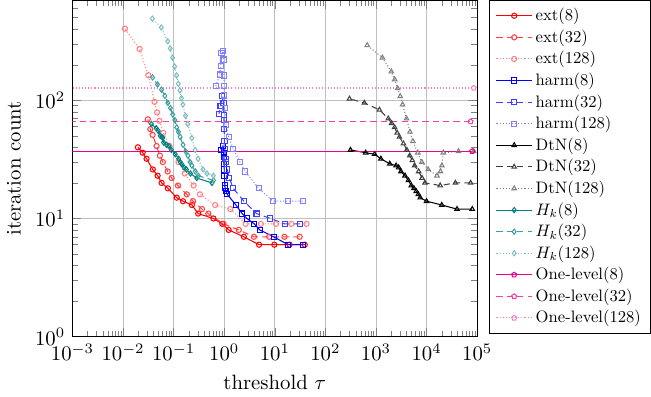}}
    \caption{
		Influence of the coarse space size (left) and threshold choice (right)
		on the iteration count in weak scaling for the imaging
		test case.
    }\label{fig:weak_scaling_vscssize_imaging}
\end{figure}

\begin{figure}[H]
    \centering
	{\includegraphics[width=0.49\textwidth]{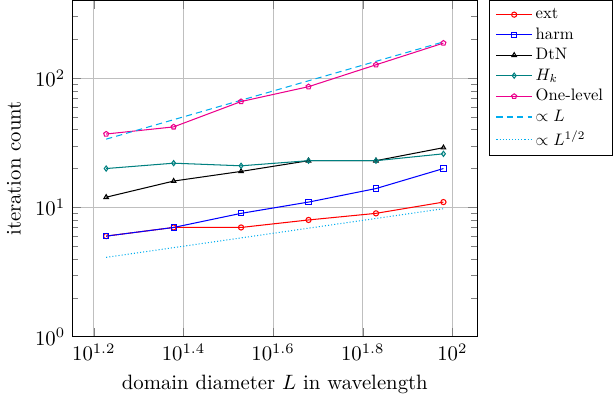}}
	{\includegraphics[width=0.49\textwidth]{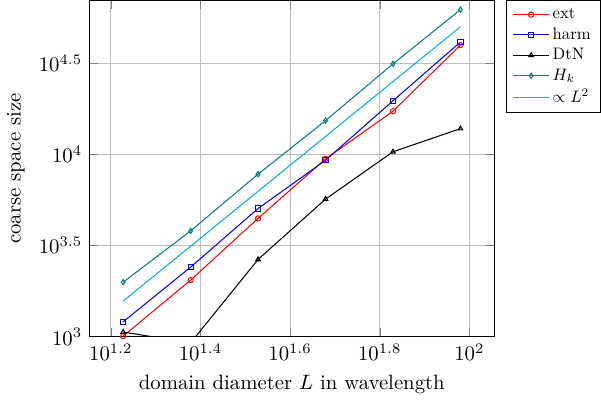}}
    \caption{
        Influence of the full domain diameter on the iteration count in weak scaling 
		for the imaging test case.
    }\label{fig:weak_scaling_vsdiameter_imaging}
\end{figure}

We then perform a weak scaling test, increasing the number of subdomains for a
fixed subdomain size.
Each subdomain is in average of size \(4.2\lambda_{0} \times 4.2\lambda_{0}\) with about
\(39000\) DOFs.
The mesh size is fixed, so the pollution effect is not taken into account.
The full results are provided in Table~\ref{tab:weak_scaling_imaging} with the
same entries as for the strong scaling test case.
The number of iterations of GMRES to reach the required tolerance as a function
of the total coarse space size and choice of threshold parameter is reported
Figure~\ref{fig:weak_scaling_vscssize_imaging}.
The number of iterations of GMRES to reach the required tolerance
as a function of the full domain diameter (measured in wavelength) is reported
Figure~\ref{fig:weak_scaling_vsdiameter_imaging}.

\paragraph{Overlap and partition of unity function}

\begin{figure}[H]
    \centering
	{\includegraphics[width=0.49\textwidth]{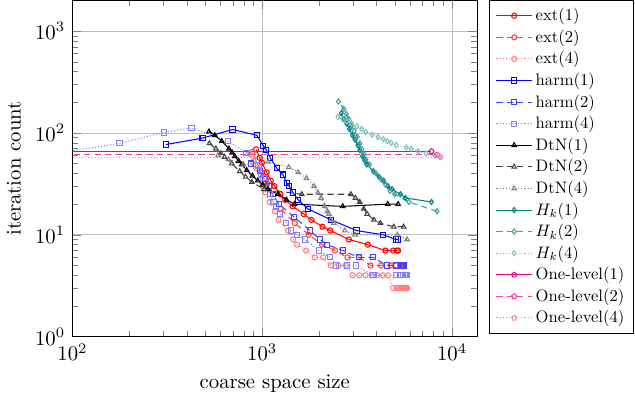}}
	{\includegraphics[width=0.49\textwidth]{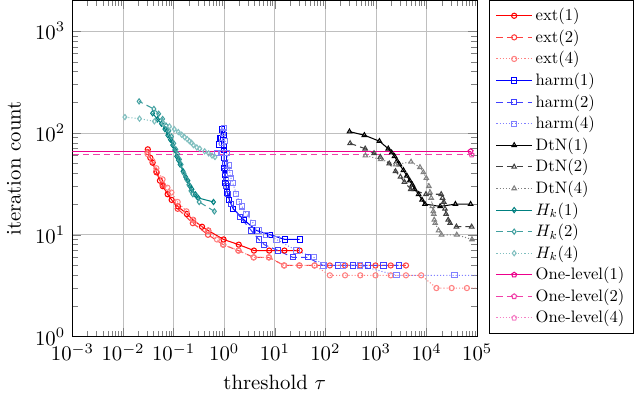}}
	{\includegraphics[width=0.49\textwidth]{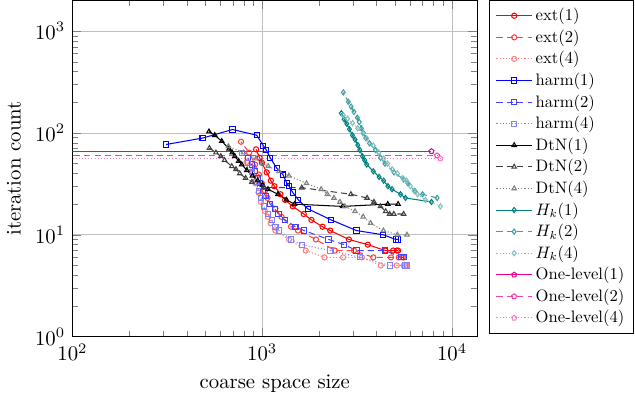}}
	{\includegraphics[width=0.49\textwidth]{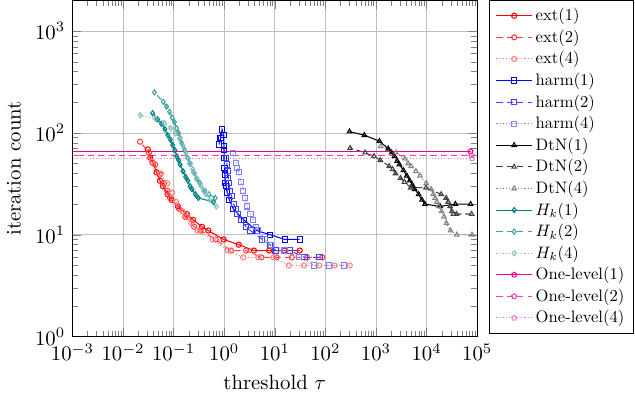}}
	{\includegraphics[width=0.49\textwidth]{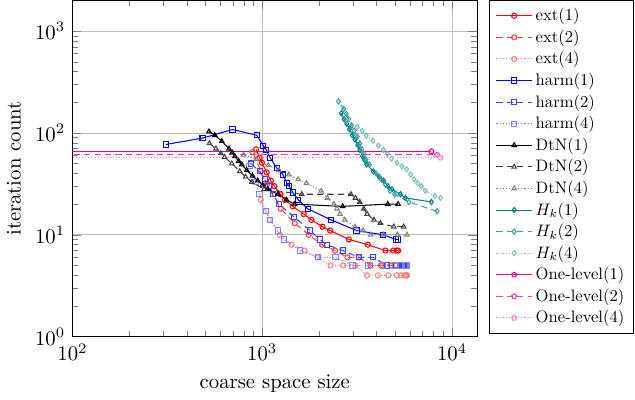}}
	{\includegraphics[width=0.49\textwidth]{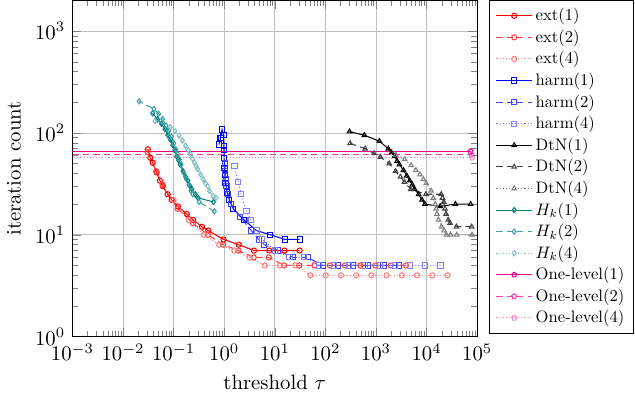}}
    \caption{
		Influence of the overlap and the partition of unity function on the
		iteration count for the imaging test case.
		Steep (top), linearly decaying (middle) and with vanishing first
		derivative on the boundary (bottom).
    }\label{fig:overlap_imaging}
\end{figure}

We turn to the study of the influence of the width of the overlap.
So far the results for this test case are obtained for minimal overlap (with a
symmetry constraint with respect to the interface).
The number of iterations of GMRES to reach the required tolerance
as a function of the total coarse space size is reported
in Figure~\ref{fig:overlap_imaging}.
The size of the overlap measured as the number of cell layers from one side of
the interface is indicated in the brackets of the legend entries.

Three different partition of unity functions, all defined as
\(\mathbb{P}_{1}\) functions hence piecewise affine, are used.
The first one is the steepest, going from \(1\) to \(0\) in the smallest number
of cells (typically two) close to the interface.
The second one is the smoothest, going from \(1\) to \(0\) on the full size of
the overlap.
The third one is an intermediate choice, decaying smoothly on the overlap but
vanishing one layer of cells before reaching the boundary, hence ensuring that
its first derivative vanishes on the boundary (which makes sense in our case
since we are using optimised boundary conditions as transmission conditions).

\clearpage
~
\vspace{0.1\textheight}

\begin{table}[H]
	\scriptsize
	\begin{center}
\begin{tabular}{ cccccc|c|ccc|ccc|ccc|ccc }
\multicolumn{6}{c|}{} & \multicolumn{1}{c|}{1lvl} & \multicolumn{3}{c|}{ext} & \multicolumn{3}{c|}{harm} & \multicolumn{3}{c|}{DtN} & \multicolumn{3}{c}{$H_k$} \\
$L$ [$\lambda_{0}$] & $N$ & $n$ & $H$ [$\lambda_{0}$] & $n_{s}$ & $n_{s}^{\partial\Omega_{s}}$ & It & It & CS & CS$_{s}$ & It & CS & CS$_{s}$ & It & CS & CS$_{s}$ & It & CS & CS$_{s}$ \\
\hline
 33.8 & 8 & 1213761 & 11.9 & 154757 & 1664 & 37 & 6 & 2073 & 259 & 6 & 2371 & 296 & 12 & 2663 & 333 & 24 & 3894 & 487 \\
 33.8 & 16 & 1213761 & 8.4 & 78134 & 1121 & 46 & 7 & 2499 & 156 & 8 & 2907 & 182 & 16 & 1807 & 113 & 23 & 5380 & 336 \\
 33.8 & 32 & 1213761 & 6.0 & 39751 & 812 & 66 & 7 & 4456 & 139 & 9 & 5048 & 158 & 19 & 2652 & 83 & 21 & 7791 & 243 \\
 33.8 & 64 & 1213761 & 4.2 & 20308 & 568 & 83 & 8 & 6460 & 101 & 11 & 7098 & 111 & 24 & 2245 & 35 & 22 & 10906 & 170 \\
 33.8 & 128 & 1213761 & 3.0 & 10487 & 413 & 112 & 9 & 8198 & 64 & 13 & 9874 & 77 & 26 & 3564 & 28 & 23 & 15849 & 124 \\
 33.8 & 256 & 1213761 & 2.1 & 5454 & 286 & 145 & 10 & 14222 & 56 & 16 & 14608 & 57 & 31 & 4878 & 19 & 24 & 22005 & 86 \\
\end{tabular}
\end{center}

    \caption{
		Strong scaling experiment for the imaging test case.
    }\label{tab:strong_scaling_imaging}
\end{table}

\clearpage
~
\vspace{0.1\textheight}

\begin{table}[H]
	\scriptsize
	\begin{center}
\begin{tabular}{ cccccc|c|ccc|ccc|ccc|ccc }
\multicolumn{6}{c|}{} & \multicolumn{1}{c|}{1lvl} & \multicolumn{3}{c|}{ext} & \multicolumn{3}{c|}{harm} & \multicolumn{3}{c|}{DtN} & \multicolumn{3}{c}{$H_k$} \\
$L$ [$\lambda_{0}$] & $N$ & $n$ & $H$ [$\lambda_{0}$] & $n_{s}$ & $n_{s}^{\partial\Omega_{s}}$ & It & It & CS & CS$_{s}$ & It & CS & CS$_{s}$ & It & CS & CS$_{s}$ & It & CS & CS$_{s}$ \\
\hline
 16.9 & 8 & 292873 & 6.0 & 38108 & 829 & 37 & 6 & 1007 & 126 & 6 & 1206 & 151 & 12 & 1060 & 132 & 20 & 1990 & 249 \\
 23.9 & 16 & 587129 & 6.0 & 38312 & 792 & 42 & 7 & 2042 & 128 & 7 & 2406 & 150 & 16 & 924 & 58 & 22 & 3804 & 238 \\
 33.8 & 32 & 1213761 & 6.0 & 39751 & 812 & 66 & 7 & 4456 & 139 & 9 & 5048 & 158 & 19 & 2652 & 83 & 21 & 7791 & 243 \\
 47.8 & 64 & 2366557 & 6.0 & 38856 & 798 & 86 & 8 & 9443 & 148 & 11 & 9288 & 145 & 23 & 5685 & 89 & 23 & 15327 & 239 \\
 67.5 & 128 & 4809133 & 6.0 & 39573 & 817 & 127 & 9 & 17250 & 135 & 14 & 19697 & 154 & 23 & 10340 & 81 & 23 & 31375 & 245 \\
 95.5 & 256 & 9520591 & 6.0 & 39238 & 811 & 187 & 11 & 39786 & 155 & 20 & 41337 & 161 & 29 & 13850 & 54 & 26 & 62267 & 243 \\
\end{tabular}
\end{center}

    \caption{
		Weak scaling experiment for the imaging test case.
    }\label{tab:weak_scaling_imaging}
\end{table}

\paragraph{Assessment and comparative conclusions for the medical imaging test case}

From these tests we can draw the following conclusions:
\begin{itemize}
	\item The two-level methods remain robust in both weak and strong scalings,
		albeit with a moderate increase in the iteration count with the number
		of subdomains, which is well-controlled by enlarging the coarse space.
	\item Harmonic and extended-harmonic coarse spaces yield the lowest
		iteration counts for a given coarse space size, with a slight advantage
		for the extended harmonic coarse space for larger number of subdomains.
	\item DtN coarse spaces remain competitive and compact, with a slightly
		worse iteration count but with the lowest size.
		We remark however some loss of robustness in the sense that the
		decrease in the number of iteration is not always monotonic with the
		size of the coarse space.
 	\item $H_k$-GenEO coarse spaces are robust but more expensive.
	\item Sensible choices for the threshold \(\tau\) are comparable to
		the previous test case, except for the DtN coarse spaces which required
		significantly larger values.
	\item Wider overlaps generally reduce iteration counts, especially for the
		harmonic and extended-harmonic coarse spaces.
	\item Smoother partition of unity functions (especially with vanishing
		derivative) yield slightly better convergence.
\end{itemize}

The imaging test case highlights the necessity of two-level methods
heterogeneous for media with overall performance similar to the simple test on a
square domain.

\subsection{Cobra cavity test case}

\paragraph{Description of the problem}

We now move to a more sophisticated test case for spectral coarse spaces: the \emph{COBRA cavity}, a three-dimensional, S-shaped waveguide originally designed for electromagnetic scattering studies by EADS Aerospatiale Matra Missiles as part of the EM-JINA 98 workshop (see~\cite{jin2015finite,liu2003scattering}). 
This geometry has since been adopted as a benchmark in domain decomposition research (e.g.\ \cite{Dolean:2015:ETC,Bonazzoli:2019:ADD}) due to its challenging features for mid- to high-frequency wave propagation.

Unlike previous test cases with simpler geometries, the COBRA cavity introduces \emph{geometric complexity} through its curvature, which can lead to \emph{wave trapping effects}. Although the wavespeed is constant, the intricate shape of the cavity poses significant numerical challenges that go beyond those of straight waveguides.

The cross-section of the cavity measures 11\,cm~$\times$~8.4\,cm. The cavity walls are modelled as sound-soft, imposing Dirichlet boundary conditions. To reduce the unbounded scattering domain to a finite computational domain, we embed the cavity in a surrounding box.
The material parameters are \(a(\boldsymbol{x}) \equiv 1\) and \(m(\boldsymbol{x}) \equiv 1\) and only the wavenumber \(k\) is varied.
All geometrical parameters are measured in terms of the wavelength \(\lambda=2\pi/k\).
The sides of the box are positioned 10\,cm from the cavity in all directions, corresponding to between 1.3 and 5.7 wavelength \(\lambda\) depending on the frequency used in the weak scaling study presented below. Robin (impedance) boundary conditions are applied on all faces of the outer box.

A normally incident plane wave excites the cavity. The problem is discretised using $\mathbb{P}_2$ Lagrange finite elements, with a resolution of 8 points per wavelength to ensure adequate accuracy in the high-frequency regime.

\begin{figure}[H]
    \centering
	%{\includegraphics[width=0.48\textwidth,clip=true,trim=7cm 0cm 7cm 0cm]{cobra_mesh.png}}
	{\includegraphics[width=0.48\textwidth,clip=true,trim=7cm 0cm 7cm 0cm]{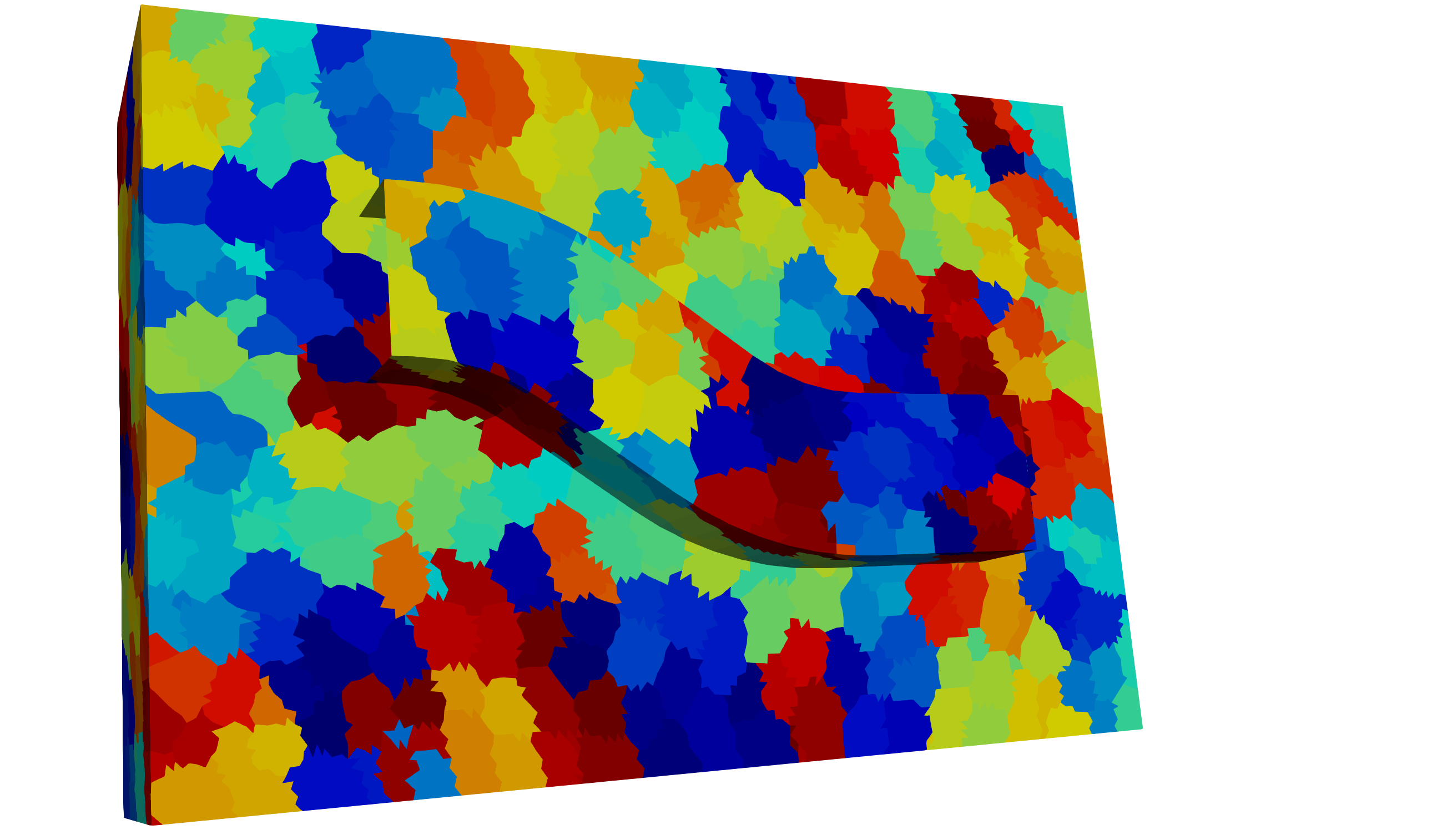}}
	{\includegraphics[width=0.48\textwidth,clip=true,trim=7cm 0cm 7cm 0cm]{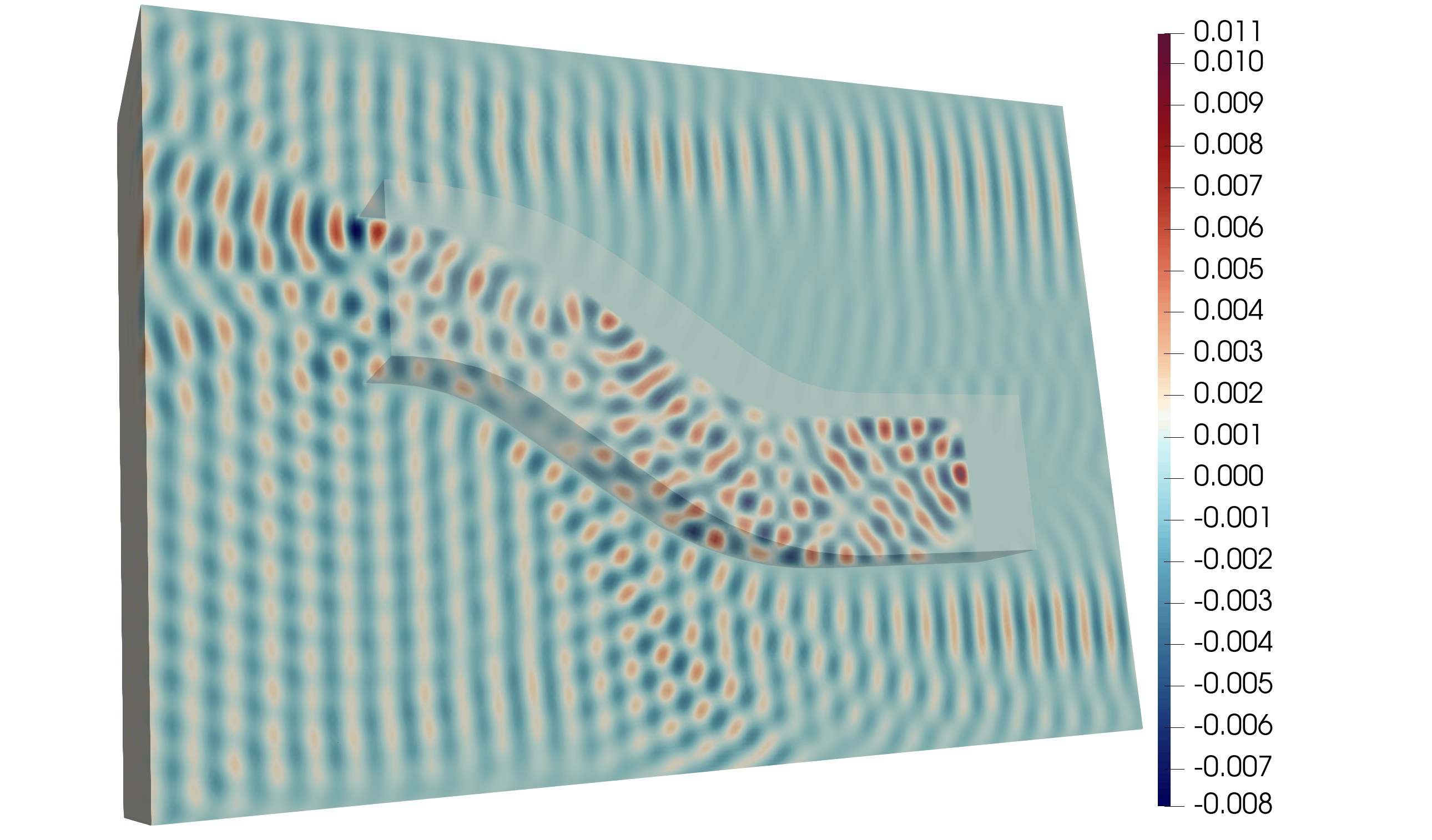}}
	%{\includegraphics[width=0.48\textwidth,clip=true,trim=7cm 0cm 7cm 0cm]{cobra_real_2.png}}
    \caption{
		Partitioning of the cobra cavity domain into 2916 subdomains (left) and real part of the total field (right) for $L = 22.9\lambda$ ($k = 360m^{-1}$). 
    }\label{fig:solutions_cobra}
\end{figure}

\paragraph{Weak scaling test}

To assess scalability in more realistic and complex 3D geometries, we perform a \emph{weak scaling test} on the COBRA cavity benchmark. In this experiment, the number of subdomains and total degrees of freedom increase proportionally with the wavenumber~$k$, while keeping the average \emph{subdomain size fixed} (approximately 1.6 wavelength \(\lambda\) in diameter and around 27{,}000 degrees of freedom per subdomain). The number of subdomains $N$ increases from 32 for $k = 80m^{-1}$ to 2916 for $k = 360m^{-1}$. Correspondingly, the diameter $L$ of the global domain in wavelength goes from \(5.1\lambda\) to \(22.9\lambda\).
As an illustration, Figure~\ref{fig:solutions_cobra} shows the partitioning of the domain as well as the real part of the solution for $L = 22.9\lambda$.

\begin{figure}[H]
    \centering
        {\includegraphics[width=0.49\textwidth]{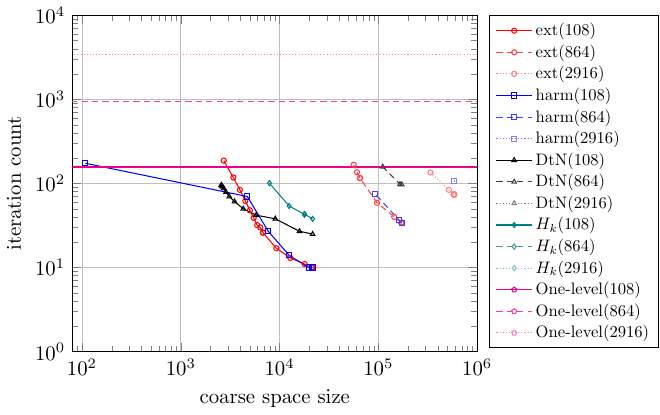}}
        {\includegraphics[width=0.49\textwidth]{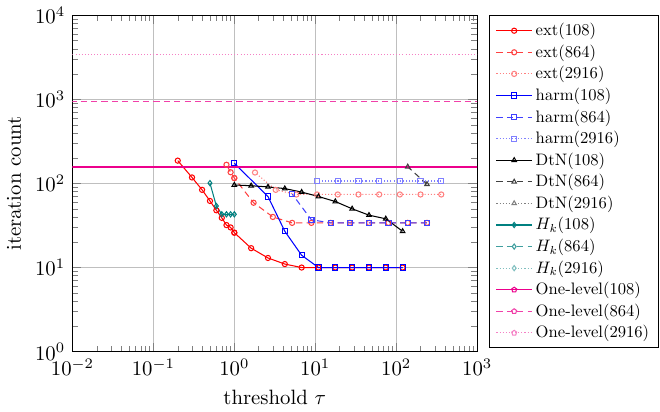}}
    \caption{
        Influence of the
		coarse space size (left)
		and threshold choice (right)
		on the iteration count
		in weak scaling 
		for the cobra cavity test case.
    }\label{fig:weak_scaling_vscssize_cobra}
\end{figure}

\begin{figure}[H]
    \centering
        {\includegraphics[width=0.49\textwidth]{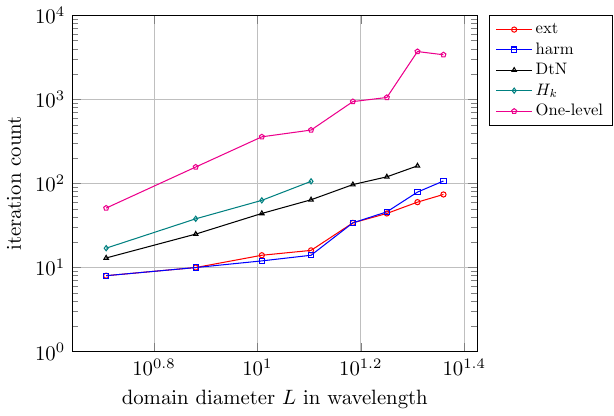}}
        {\includegraphics[width=0.49\textwidth]{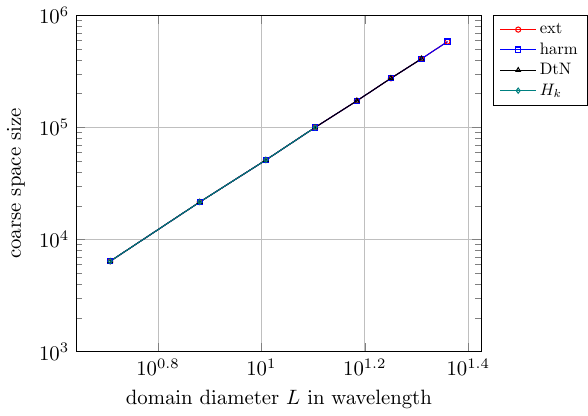}}
    \caption{
        Influence of the full domain diameter on the iteration count in weak scaling 
		for the cobra cavity test case.
    }\label{fig:weak_scaling_vsdiameter_cobra}
\end{figure}

As shown in Table~\ref{tab:weak_scaling_cobra}, the domain diameter measured in wavelengths increases from $L = 5.1\lambda$ (for $N = 32$) up to $L = 22.9\lambda$ (for $N = 2916$). This setup preserves local resolution while increasing the global problem size, making it ideal to evaluate the robustness of coarse spaces and the scalability of the two-level ORAS preconditioner in the mid- to high-frequency regime.
In addition to that, we also report the number of GMRES iterations to reach the required tolerance as a function of the total coarse space size and choice of threshold parameter in Figure~\ref{fig:weak_scaling_vscssize_cobra}.
The number of subdomains is provided in brackets in the legend entries.
The number of iterations of GMRES to reach the required tolerance as a function of the full domain diameter (measured in wavelength) is reported in Figure~\ref{fig:weak_scaling_vsdiameter_cobra}.

\begin{table}[H]
    \scriptsize
    \begin{center}
\begin{tabular}{ cccccc|c|ccc|ccc|ccc|ccc }
\multicolumn{6}{c|}{} & \multicolumn{1}{c|}{1lvl} & \multicolumn{3}{c|}{ext} & \multicolumn{3}{c|}{harm} & \multicolumn{3}{c|}{DtN} & \multicolumn{3}{c}{$H_k$} \\
$L$ [$\lambda_{}$] & $N$ & $n$ & $H$ [$\lambda_{}$] & $n_{s}$ & $n_{s}^{\partial\Omega_{s}}$ & It & It & CS & CS$_{s}$ & It & CS & CS$_{s}$ & It & CS & CS$_{s}$ & It & CS & CS$_{s}$ \\
\hline
5.1  &   32 & 473004   & 1.6 & 22336 & 2986 & 51   &  8 &   6400 & 200 &  8 &   6400 & 200 & 13   &   6400 & 200 &   17 &   6400 & 200 \\
7.6  &  108 & 1761164  & 1.6 & 25895 & 3824 & 157  & 10 &  21600 & 200 & 10 &  21600 & 200 & 25   &  21600 & 200 &   38 &  21600 & 200 \\
10.2 &  256 & 4140366  & 1.6 & 26525 & 4144 & 359  & 14 &  51200 & 200 & 12 &  51200 & 200 & 44   &  51200 & 200 &   63 &  51200 & 200 \\
12.7 &  500 & 8431281  & 1.6 & 27980 & 4468 & 432  & 16 & 100000 & 200 & 14 & 100000 & 200 & 64   & 100000 & 200 & 106 & 100000 & 200 \\
15.3 &  864 & 13927097 & 1.6 & 27839 & 5035 & 942  & 34 & 172800 & 200 & 34 & 172800 & 200 & 97   & 172800 & 200 & >200 & 172800 & 200 \\
17.8 & 1372 & 18879654 & 1.6 & 24589 & 4676 & 1055 & 44 & 274400 & 200 & 46 & 274400 & 200 & 120  & 274400 & 200 & >200 & 274400 & 200 \\
20.4 & 2048 & 32848020 & 1.6 & 28143 & 5210 & 3711 & 60 & 409600 & 200 & 79 & 409600 & 200 & 162  & 409600 & 200 & >200 & 409600 & 200 \\
22.9 & 2916 & 44520439 & 1.6 & 27138 & 5121 & 3398 & 74 & 583200 & 200 & 107 & 583200 & 200 & >200 & 583200 & 200 & >200 & 583200 & 200 \\
\end{tabular}
\end{center}

    \caption{
		Weak scaling experiment for the cobra cavity test case.
    }\label{tab:weak_scaling_cobra}
\end{table}

\paragraph{Assessment and comparative conclusions for the cobra cavity test case}

\begin{itemize}
	\item GMRES iteration counts show that all two-level methods significantly
		outperform the one-level baseline, especially as the global domain size
		$L$ increases.
	\item Harmonic and extended harmonic coarse spaces consistently achieve the
		lowest iteration counts, with only mild growth across scales.
	\item The DtN coarse space performs competitively at moderate scales but
		exhibits a slight increase in iteration count at higher $L$, and
		eventually breaks down for the largest problem.
	\item \(H_{k}\)-GenEO also performs competitively at moderate scales but breaks
		down beyond a certain coarse space size and for larger domains.
\end{itemize}

\subsection{The GO\_3D\_OBS Crustal Geomodel}

\paragraph{Description of the problem}

The \texttt{GO\_3D\_OBS} crustal geomodel is a high-resolution, three-dimensional synthetic model specifically designed to evaluate seismic imaging techniques, particularly for deep crustal and subduction zone exploration~\cite{Gorszczyk_2021_GNT}. It represents a continental margin at a regional scale, incorporating realistic geological heterogeneities and sharp velocity contrasts that are characteristic of tectonic plate boundaries.

The model features a wide range of acoustic wavespeeds, \(c(\boldsymbol{x})\) varies from $1500\,\mathrm{ms^{-1}}$ in the near-surface sedimentary layers and water column up to $8639\,\mathrm{ms^{-1}}$ in the lower crust and upper mantle. This range captures the physical complexity required to simulate realistic wave propagation in crustal-scale full-waveform inversion (FWI) settings.

Importantly, the geomodel includes the essential structural features of subduction zones, such as a dipping slab, accretionary prism, crust-mantle transitions, and bathymetry. It is inspired by a real-world FWI case study performed in the eastern Nankai Trough~\cite{Gorszczyk_2017_TRW}, and is designed to stress-test forward and inverse modelling tools in the presence of complex multiscale heterogeneity.

\begin{figure}[H]
\begin{center}
\includegraphics[width=0.9\textwidth,clip=true,trim=2cm 5cm 5cm 8cm]{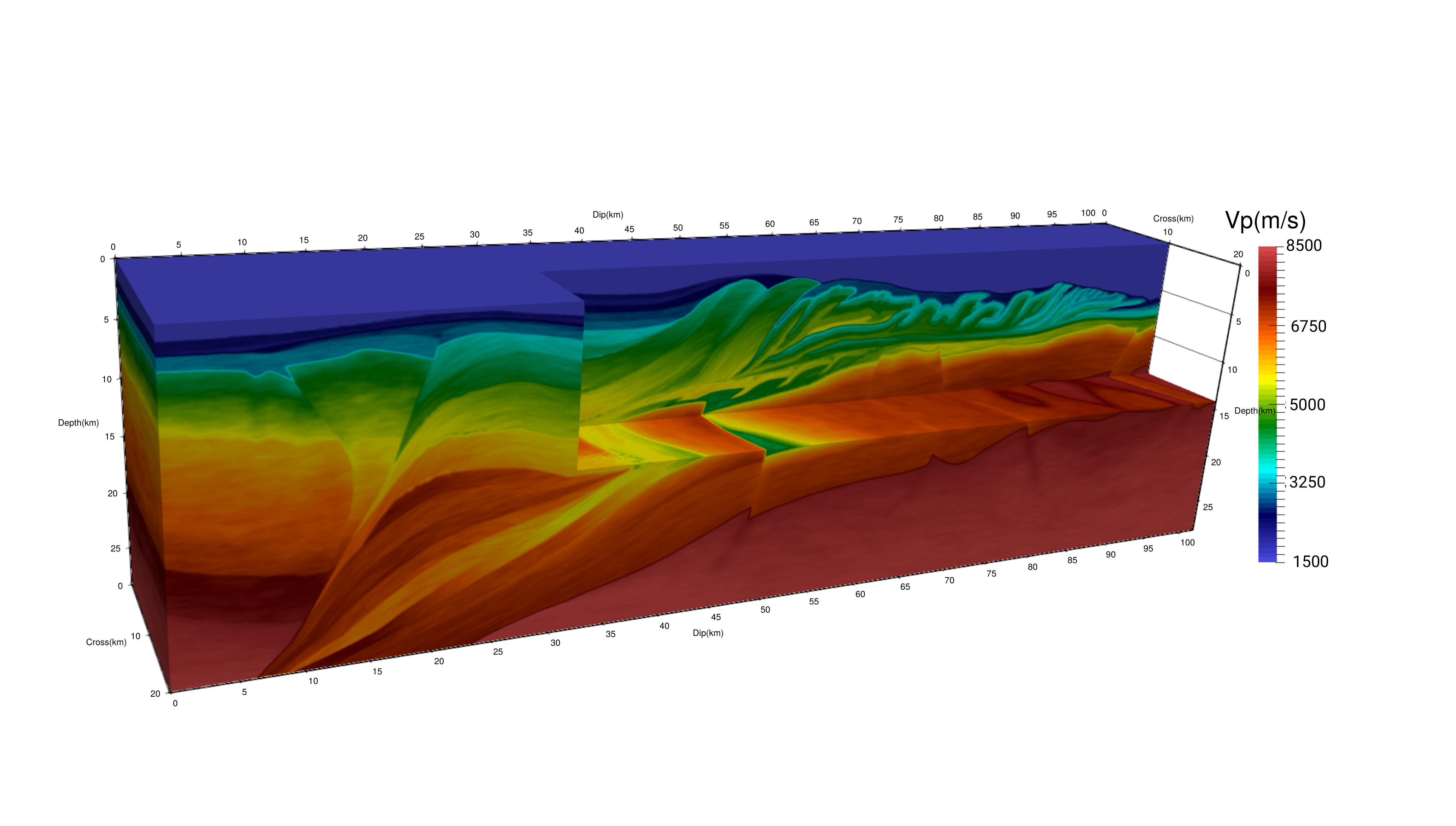}
\caption{Target of the regional GO\_3D\_OBS model representing the crust of a subduction zone \cite{Gorszczyk_2021_GNT}.}
\label{fig:go3dobs_model}
\end{center}
\end{figure}

For the simulations considered in this work, we extract a target subregion of the full model, covering $20\,\mathrm{km} \times 102\,\mathrm{km} \times 28.3\,\mathrm{km}$ (Figure~\ref{fig:go3dobs_model}).
The material coefficients of our model for this test case are set to \(a \equiv 1\) and \(m(\boldsymbol{x}) = c^{-2}(\boldsymbol{x})\).
The geometrical parameters are measured according to the minimal wavelength \(\lambda_{\min} = 2\pi\min_{\boldsymbol{x}}(c(\boldsymbol{x}))/\omega\).
The domain is discretised with an unstructured tetrahedral mesh, adapted to the local wavelength to ensure accurate resolution of the wavefield, with denser refinement near the bathymetry and across high-contrast interfaces.

\begin{figure}[H]
\begin{center}
\includegraphics[width=0.9\textwidth]{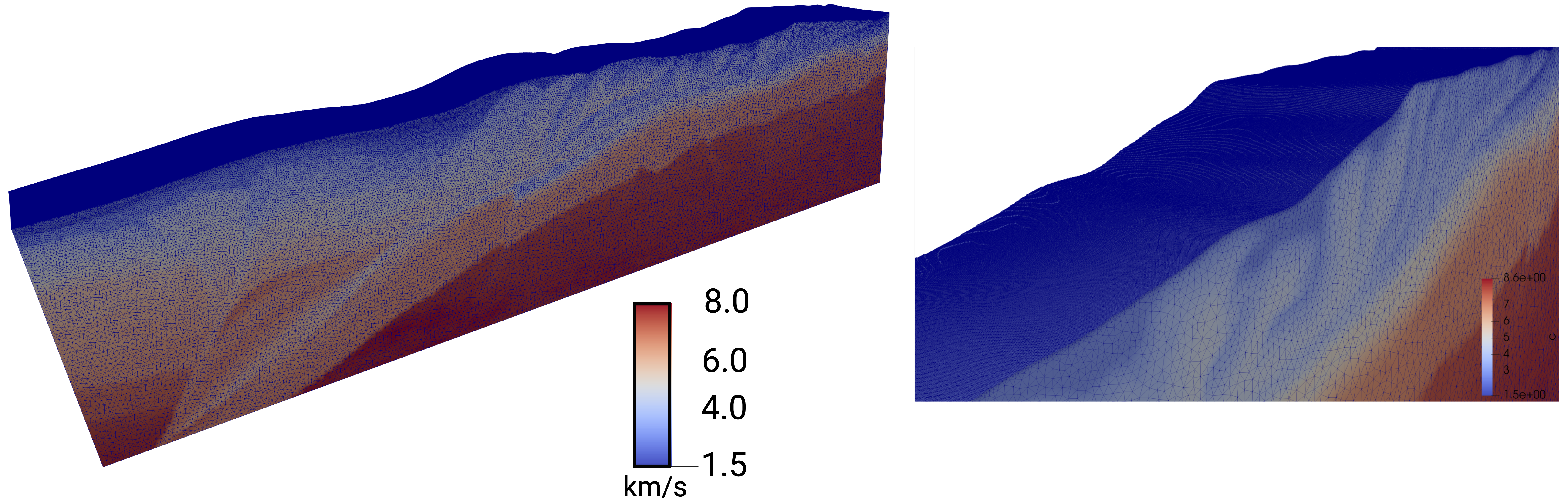}
\caption{Tetrahedral mesh of the GO\_3D\_OBS model adapted to the local wavelength. Note that the mesh in the water layer is not shown in the Figure. Instead, the conformal meshing of the seabed is highlighted. The right panel shows a zoom.}
\label{fig:go3dobs_mesh}
\end{center}
\end{figure}

For this simulation, we use an unstructured mesh that complies with the bathymetry. Below the bathymetry, the size of the elements is set according to the local wavelength such that we have approximately 4 points per wavelength for the reference angular frequency $\omega = 3.75\pi\,\mathrm{rad\,s^{-1}}$ (Figure~\ref{fig:go3dobs_mesh}).

\paragraph{Strong scaling test}

To further assess scalability in heterogeneous 3D media, we perform a \emph{strong scaling test} on the crustal-geomodel benchmark.
In this experiment, we consider two partitions with \(1024\) and \(2048\) subdomains.

\begin{figure}[H]
    \centering
        {\includegraphics[width=0.49\textwidth]{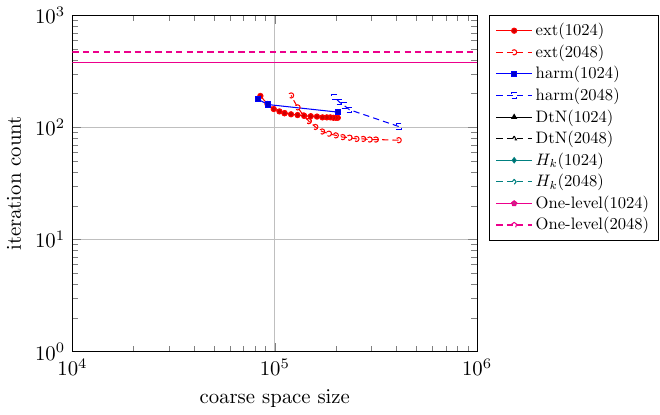}}
        {\includegraphics[width=0.49\textwidth]{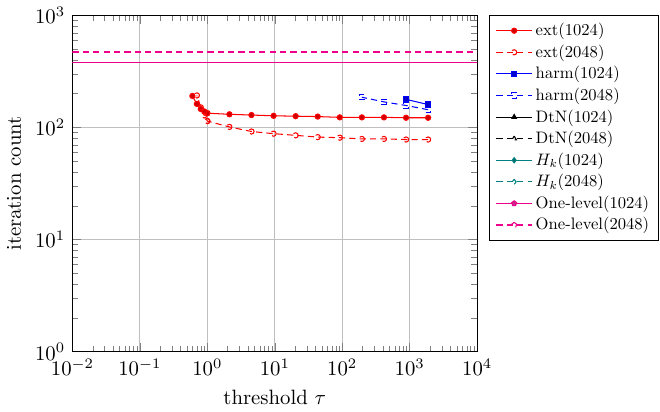}}
    \caption{
        Influence of the
		coarse space size (left)
		and threshold choice (right)
		on the iteration count
		in a strong scaling experiment
		for the GO\_3D\_OBS test case.
    }\label{fig:strong_scaling_vscssize_go3dobs}
\end{figure}

The results are reported in Table~\ref{tab:strong_scaling_go3dobs} and we also
report the number of GMRES iterations to reach the required tolerance as a
function of the total coarse space size and choice of threshold parameter in
Figure~\ref{fig:strong_scaling_vscssize_go3dobs}.
The number of subdomains is provided in brackets in the legend entries.

\begin{table}[H]
	\scriptsize
	\begin{center}
\begin{tabular}{ cccccc|c|ccc|ccc|ccc|ccc }
\multicolumn{6}{c|}{} & \multicolumn{1}{c|}{1lvl} & \multicolumn{3}{c|}{ext} & \multicolumn{3}{c|}{harm} & \multicolumn{3}{c|}{DtN} & \multicolumn{3}{c}{$H_k$} \\
$L$ [$\lambda_{\min}$] & $N$ & $n$ & $H$ [$\lambda_{\min}$] & $n_{s}$ & $n_{s}^{\partial\Omega_{s}}$ & It & It & CS & CS$_{s}$ & It & CS & CS$_{s}$ & It & CS & CS$_{s}$ & It & CS & CS$_{s}$ \\
\hline
 48.3 & 1024 & 22297073 & 4.8 & 35185 & 5416 & 380 & 122 & 204800 & 200 & 137 & 204800 & 200 & >200 & 204800 & 200 & >200 & 204800 & 200 \\
 48.3 & 2048 & 22297073 & 3.8 & 19864 & 3712 & 471 & 77 & 409600 & 200 & 102 & 409600 & 200 & >200 & 409600 & 200 & >200 & 409600 & 200 \\
\end{tabular}
\end{center}
    \caption{
		Strong scaling experiment for the GO\_3D\_OBS test case.
    }\label{tab:strong_scaling_go3dobs}
\end{table}

\paragraph{Assessment and comparative conclusions for the crustal-geomodel test case}

\begin{itemize}
	\item Harmonic and extended-harmonic coarse spaces are able to provide a
		reduction in the iteration count compared to the one-level method. 
		The smallest number of iterations is obtained for the extended-harmonic
		one, with more compact coarse spaces in the case with the largest
		number of subdomains.
	\item In contrast, both DtN and \(H_{k}\)-GenEO coarse spaces do not provide
		convergence for this challenging heterogeneous 3D test case.
	\item The harmonic coarse space required significantly larger values for
		the eigenvalue threshold parameter.
\end{itemize}

\section{Conclusions}
\label{sec:Conclusions}

This study presents a comprehensive evaluation of spectral coarse spaces in two-level overlapping Schwarz methods for solving high-frequency Helmholtz problems in both homogeneous and heterogeneous media. The analysis spans a wide range of numerical settings, from idealised 2D test cases to large-scale 3D benchmarks such as the COBRA cavity and the GO\_3D\_OBS crustal geomodel.

A comparative overview of the key characteristics of each spectral coarse space is provided in Table~\ref{tab:cs_summary}.
Our main findings can be summarised as follows:

\begin{itemize}
	\item \textbf{Two-level domain decomposition is essential} for robust and
		scalable solvers at mid-to-high frequencies. All spectral coarse spaces
		studied outperform the one-level ORAS baseline by a wide margin.
  
	\item \textbf{Extended harmonic coarse spaces} offer the best trade-off
		between solver efficiency and coarse space size. They consistently
		deliver low GMRES iteration counts, even under weak scaling and in
		high-frequency 3D heterogeneous settings. We point out however that the
		eigenproblem is posed in an enlarged domain.
  
	\item \textbf{Harmonic coarse spaces} have a performance comparable to that of
		extended-harmonic coarse spaces.
		In some occurrences, in particular with a large number of subdomains,
		they require slightly larger coarse space sizes or a few more
		iterations compared to their extended counterpart.
  
	\item \textbf{DtN coarse spaces} are compact and competitive for moderate
		problem sizes but become less effective for larger domains and higher
		frequency, and breakdown can be observed in 3D.
		We remark a strong sensitivity with respect to the eigenvalue threshold
		parameter \(\tau\).
  
	\item \textbf{{\boldmath\(H_{k}\)}-GenEO} exhibits strong robustness in 2D but incurs
		significantly higher computational cost due to larger coarse spaces,
		and breakdown can be observed in 3D. It remains an attractive choice
		when memory constraints are less critical.

	\item Additional implementation aspects---such as overlap width, partition
		of unity smoothness, and eigenvalue threshold tuning---affect
		performance and should be carefully optimised in practice.
\end{itemize}

Overall, the results validate the use of physics-informed spectral coarse spaces---especially harmonic and extended harmonic types---as effective strategies for achieving scalability and wavenumber robustness in Helmholtz solvers across a wide spectrum of applications.

A comparison of these domain decomposition preconditioners in terms of setup
and GMRES run times in addition to iteration count and coarse space sizes has
been left outside of the scope of this study but is of practical importance and
could be a perspective for future work.

\begin{table}[H]
	\centering
	\footnotesize
	\begin{tabu}{c|c|c|c|c|c|c|c|c}
		\textbf{Problem} & $d$ & \textbf{Medium} & \textbf{B.C.} & \textbf{Scaling} & \textbf{Extended} & \textbf{Harmonic} & \textbf{DtN} & \textbf{$H_{k}$-GenEO} \\
		\hline
		\multirow{2}{*}{Square}           & 2D & homogeneous   & Robin & strong & \tick & \tick & \tick\tick & \tick \\
		                                  & 2D & heterogeneous & Robin & strong & \tick & \tick & \tick\tick & \tick \\
		\hline
		\multirow{2}{*}{Medical imaging}  & 2D & heterogeneous & PML   & strong & \tick\tick & \tick\tick & \tick & \tick \\
		                                  & 2D & heterogeneous & PML   & weak   & \tick\tick & \tick\tick & \tick & \tick \\
		\hline
		Cobra cavity                      & 3D & homogeneous   & Robin & weak   & \tick\tick & \tick\tick & \tick/\cross & \tick/\cross \\
		\hline
		Crustal geomodel                  & 3D & heterogeneous & PML   & strong & \tick\tick & \tick\tick & \cross & \cross
	\end{tabu}
	\caption{
		An overview of which coarse spaces perform well in the different
		problem scenarios tests.
		A \tick~indicates that the method performs well, with \tick\tick~indicating
		this method was most favourable in a particular instance.
		A \cross~indicates that a method provided relatively little to no gain over
		the one-level method.
	}
	\label{tab:cs_summary}
\end{table}

\printbibliography[title={References}]

\end{document}